\documentclass[reqno]{amsart}

\usepackage[latin1]{inputenc}
\usepackage[T1]{fontenc}                  
\usepackage[english]{babel}
\usepackage{amsmath,amssymb,amsbsy,amsthm,amsfonts}
\usepackage{color}
\usepackage{epsfig}
\usepackage{pgf}
\usepackage{subfigure}

\newtheorem{definition}{Definition}[section]
\newtheorem{remark}[definition]{Remark}

\newtheorem{theorem}[definition]{Theorem}
\newtheorem{proposition}[definition]{Proposition}

\newtheorem{lemma}[definition]{Lemma}

\usepackage[colorlinks=true]{hyperref}
\usepackage{theoremref}

\usepackage{accents}

\usepackage{wasysym}

\newcommand{\eps}{\varepsilon}
\newcommand{\R}{\mathbb{R}}

\def\tilde{\widetilde}

\usepackage{enumerate}
\usepackage{xcolor}

\begin{document}
\title[Stationary solutions to the PNP equations with steric effects] 
      {Stationary solutions to the Poisson-Nernst-Planck equations with steric effects}  
\author[L.-C. Hung and X. Liao]{Li-Chang Hung$^{\ast}$ and Mach Nguyet Minh$^{\natural}$}

\email{lichang.hung@gmail.com; minh.mach@helsinki.fi}
\thanks{$^{\ast}$Department of Mathematics, National Taiwan University Taiwan}
\thanks{$^{\natural}$Department of Mathematics and Statistics, University of Helsinki, Finland}

\keywords{Poisson-Nernst-Planck equations, steric effects, steady state, differential algebraic equations, unique solution, bifurcation, semilinear Poisson equation, homogeneous Neumann boundary conditions. }


\begin{abstract}
Ion transport, the movement of ions across a cellular membrane, plays a crucial role in a wide variety of biological processes and can be described by the Poisson-Nernst-Planck equations with steric effects (PNP-steric equations). In this paper, we shall show that under homogeneous Neumann boundary conditions, the steady-state PNP-steric equations are equivalent to a system of differential algebraic equations (DAEs).  Analyzing this system of DAEs inspires us to propose an assumption on coupling constants, the so-called \textbf{(H1)} which will be introduced in Section~\ref{Sec:model}, such that if \textbf{(H1)} holds true, the steady-state PNP-steric equations admit a unique stationary  $C^2$ solution. Moreover, we shall point out the occurrence of bifurcation when \textbf{(H1)} is violated, which may relate to the opening and closing of the ion channels. When \textbf{(H1)} fails, we also suggest a simple criterion to check whether the system of DAE equations admits unique  monotone $C^2$ solutions; or unique  monotone piecewise $C^2$ solutions with vertical tangents; or triple piecewise $C^2$ solutions.  To the best of the authors' knowledge, this is the first time such DAE approach has been utilized to  obtain a complete investigation for the steady-state PNP-steric equations of two counter-charged ion species.
\end{abstract}

\maketitle

\section{Introduction}\label{sec: intro}

The Poisson-Nernst-Planck (PNP) equations have been used to describe the diffusion of charged particles under the influence of an electric field since the past century, and have a wide range of applications: from electrochemistry \cite{rubinstein1990electro,durand1996new,barthel1998physical,bazant2004diffuse,fawcett2004liquids,bazant2005current} to the semiconductor devices \cite{roulston1990bipolar,warner2001microelectronics,streetman2005solid,selberherr2012analysis,jerome2012analysis}. 

In biophysics, the PNP equations were  suggested as the basic continuum model for simulating the movements of ions across the cellular membrane through open ion channels
\cite{chen1993charges,eisenberg1996computing,eisenberg1998ionic}.  
Ion channels are pore-forming proteins located in the cellular membrane that have the ability to open and close in response to chemical or mechanical signals. Ion channels can be also called passageways, since they usually allow only a single type of ion to pass through them. Therefore, ion channels play an essential role in cell sustaining and control many important physiological processes such as nerve and muscle excitation, cell volume and blood pressure regulation, cell proliferation, hormone secretion, fertilisation, learning and memory, programming cell death \cite{hacker2009pharmacology}. Also, {detailed knowledge} of ion channels {is} {very useful} for new drug design and efficient gene therapy \cite{niemeyer2001ion}.
The continuum PNP equations were derived from a Langevin model of ionic motion \cite{shuss2001derivation,nadler2004ionic} and can be considered as the most simplified and successful model for ion flow through membrane channel compared to others such as {ab initio} molecular dynamics and classical molecular dynamics, since the continuum PNP equations are able to yield {good predictions} of ion channel transport at a relatively small computational cost \cite{zheng2011second}. 

Over the last decades, a wide range of computational algorithms, including finite difference, finite element and finite volume methods, have been proposed for the numerical solutions of the PNP equations. A fully self-consistent numerical solution of the PNP equations for a cylindrical channel in 3D was first studied in \cite{barcilon1992ion,chen1992constant}. The authors in \cite{kurnikova1999lattice} then developed a lattice relaxation algorithm in combination with the finite difference method for solving the PNP equations through arbitrary 3D volumes and also gave an application {to} the 3D realistic geometry of the Gramicidin A channel. {In \cite{hollerbach2001two}, the authors also considered} the Gramicidin A channel and made use of the spectral element method, which is a {particular version of} finite element method, that allows to employ { more physically meaningful} boundary conditions. Convergence was substantially improved in \cite{mathur2009multigrid} through the use of a Newton-Raphson iteration procedure coupled to an algebraic multigrid method and an unstructured cell-centered finite volume method discretization. The first second-order convergent numerical scheme for solving PNP equations for realistic ion channels was introduced in \cite{zheng2011second}.


In contrast to the {numerous works} on the numerical algorithms, {there are a few} results that related to the mathematical aspects of the PNP equations {in the literature}.
Due to the strongly coupled equations, complexity of the irregular geometry and presence of geometric singularities, full scale mathematical analysis of the PNP equations such as existence, uniqueness, asymptotic behavior, stability as well as analytic formula of the solutions, under realistic biological setting is highly challenging and yet to be achieved. The existence and stability of the solutions of the steady-state PNP equations for electron flows in semiconductors, on the other hand, were established in \cite{jerome1985consistency}, and the existence and long time behavior of the unsteady PNP equations were studied in \cite{biler1994debye}.

Although the PNP equations {give good predictions} of experimental measurements of ion transport problems, the continuum PNP model itself still contains many limitations. Indeed, based on a mean-field approximation of ions, the continuum PNP model treats ions as continuous charge densities. As a consequence, the finite volume effect of ion particles and non-electrostatic interactions of ion species are neglected in the PNP theory. Moreover, the PNP model also {lacks of} the description of ionic dielectric boundary effects. To address these drawbacks, many modified PNP models have been proposed in the literature \cite{gillespie2002coupling,corry2003dielectric,di2004specific,kilic2007steric,jung2009computational,li2009continuum,li2010solutions,bazant2011double,lu2011poisson,horng2012pnp,lin2013poisson,lin2014new,liu2014poisson}. The reader is referred to \cite{iglivc2015nanostructures} for an overview and to \cite{chung2002recent,roux2004theoretical,coalson2005poisson} for discussions about advantages and limitations of these modified PNP models. Among these models, we shall follow the results in \cite{lin2014new} which focused on the ion-size effects (aka. steric effects) caused by finite size ions crowded in a narrow channel. The mathematical model for the PNP equations with steric effects (PNP-steric equations) proposed in \cite{lin2014new} is actually a simplified version of the one in \cite{Eisenberg-Liu-11-hard-sphere-repulsion-ionic,Eisenberg-Liu-10-Energy-variational-ions}. Applying the idea from liquid state theory, the authors {in} \cite{Eisenberg-Liu-11-hard-sphere-repulsion-ionic,Eisenberg-Liu-10-Energy-variational-ions} modified the continuum PNP equations by adding the repulsive term of the Lennard-Jones (LJ) potential to the energy functional of the PNP equations. The LJ potential is a well-known mathematical model for describing the interaction between a pair of ions and often used as an approximate model of the van der Waals force \cite{parsegian2005van}. However, since the LJ potential is singular at the origin, the modified PNP equations in \cite{Eisenberg-Liu-11-hard-sphere-repulsion-ionic,Eisenberg-Liu-10-Energy-variational-ions} become a complicated system of differential-integral equations with singular integrals which allows no theoretical result. Besides, numerical solutions may become inaccurate because of the effect of high Fourier frequencies \cite{Eisenberg-Liu-11-hard-sphere-repulsion-ionic}. By approximating the LJ potential with band-limited functions, the authors in  \cite{lin2014new} obtained the PNP-steric equations which merely contain nonlinear differential terms (with coupling constants) instead of singular integrals. \cite{lin2014new} also provided the stability and instability conditions for the PNP-steric equations in 1D, with two species and with zero Dirichlet boundary conditions. Numerical efficiency of the PNP-steric equations was shown in \cite{horng2012pnp}.  

In this paper, we focus on the steady-state solutions of the PNP-steric equations. The steady-state solutions are obtained by setting all the time derivatives in the PNP-steric equations to zero, and can be considered as the first step in order to understand the asymptotic behavior of the system from a physical point of view. The existence of multiple steady-state solutions in 1D with Robin boundary conditions for three and four species under the assumptions that the first two species have the same coupling constants and opposite sign of valences was investigated in \cite{lin2015multiple}. Our work then {relaxes the assumptions of \cite{lin2015multiple}} on the coupling constants and the valences, and points out when the PNP-steric equations for two species of anions and cations under the homogeneous Neumann boundary conditions admit unique or multiple stationary solutions.

The rest of the paper is organized as follows. In Section~\ref{Sec:model}, we shall explain the derivation of the PNP-steric equations from the continuum PNP equations and sumarize our main results. Section~\ref{sec: Two species equations} is devoted to verify the equivalence of the steady-state PNP-steric equations and a corresponding system of DAEs. The existence and some basic properties of solutions to this DAEs system are also considered in this section. In Section~\ref{sec:H1}, we assume  {\bf (H1)} and investigate the uniqueness and $C^2$-smoothness of solutions  to the DAEs system; some more properties of solutions to the DAEs system; as well as the uniqueness and $C^2$-smoothness of solutions to the steady-state PNP-steric equations under homogeneous Neumann boundary conditions in one dimensional space. Section~\ref{sec: Triple solutions} studies the bifurcation of the solutions to the DAEs system as {\bf (H1)} is violated and analyzes when the unique and the triple piecewise $C^2$ solutions occur. We end the paper by proving some auxiliary lemmas in Section~\ref{sec:appendix}.

\section{Mathematical model and main results}\label{Sec:model}


In this section, we shall  {summarize} the results in \cite{lin2014new} {which explain} the derivation of the PNP-steric equations from the continuum PNP equations. The continuum PNP equations are {formed} by coupling the Nernst-Planck equations, which {describe} the rate of change of the concentration of each ion species due to the concentration flux of ion {diffusitivity} and electrostatic force,
\begin{equation*}
\begin{cases}
\displaystyle\frac{\partial c_i}{\partial t} = - \nabla \cdot J_i, &i=1,\dots,N, \\
J_i = - D_i\left(\nabla c_i + \displaystyle\frac{z_i\,e}{k_B\,T}\,c_i\,\nabla \phi \right),  &i=1,\dots,N,
\end{cases}
\end{equation*}
%
%
with the electrostatic Poisson equation
\[
-\nabla \cdot (\eps \nabla \phi) = \rho_0 + \sum_{i=1}^N {z_iec_i}.
\]
Here, $N$ denotes the number of ion species; $c_i$, $J_i$, $D_i$, and $z_i$  are respectively the concentration, concentration flux, diffusion constant and valence of the {$i^{\rm th}$} ion species. The electrostatic potential is denoted by $\phi$, whilst $k_B$ is the Boltzmann constant, $T$ is the absolute temperature, $e$ is the elementary charge, $\eps$ is the dielectric constant and $\rho_0$ is the permanent (fixed) charge.

The PNP-steric equations are obtained by {adding an} approximation of the repulsive term of the LJ potential {to} the concentration flux: 
\[
J^{steric}_i  = - D_i\,\nabla c_i - \frac{D_i\,c_i}{k_B\,T}\,z_i\,e\,\nabla \phi - \frac{D_i\,c_i}{k_B\,T}\,S_\sigma\,\sum_{j=1}^N \epsilon_{ij}\,(a_i+a_j)^{{12}} \nabla c_j, \qquad   i=1,\dots,N,
\]
where $S_\sigma:= \displaystyle\frac{\omega_d}{12-d}\sigma^{d-12}$, $d$ is the dimension of the {considered} Euclidean space, $\omega_d$ is the surface area of the $d$-dimensional unit ball, and $\sigma$ is the small parameter used in the spatially band-limited function to define the radius of the truncation frequency range. When $\sigma$ tends to zero, the approximate LJ potential tends to the original LJ potential. As a consequence, the total energy functional in the PNP-steric equations tends to the one in the papers \cite{Eisenberg-Liu-11-hard-sphere-repulsion-ionic} and \cite{Eisenberg-Liu-10-Energy-variational-ions}. The radius of the {$i^{\rm th}$} ion species {is now taken into account} and denoted by $a_i$, {whilst} $\epsilon_{ij}$ is an appropriately chosen energy constant which comes from the repulsive part of the LJ potential to describe the hard sphere repulsion of ions. For the notation convenience, we have assumed that $\epsilon_{ij}=\epsilon_{ji}$.

We consider a bounded domain {$\Omega \subset \mathbb{R^d}, (d\ge 1)$} with smooth boundary and the case of two counter-charged ion species (i.e. $N=2$). {The indices $i=1,2$ are now changed} to $i=n,p$ to indicate the anionic and cationic species, respectively. Denote by $\tilde{S}_\sigma:=\displaystyle\frac{S_\sigma}{k_B\,T}$, {$g_{nn}:=\epsilon_{nn}(2\,a_n)^{12}$, $g_{np}:=\epsilon_{np}(a_n+a_p)^{12}$, and $g_{pp}:=\epsilon_{pp}(2\,a_p)^{12}$}, the PNP-steric equations become
\begin{equation*}
\begin{cases}
\displaystyle\frac{\partial c_n}{\partial t} = D_n \left[\nabla \cdot \left(\nabla c_n + \frac{z_n\,e}{k_B\,T}\,c_n\,\nabla \phi \right) + \tilde{S_\sigma} \nabla \cdot \left(g_{nn}\,c_n\,\nabla c_n + g_{np}\,c_n\,\nabla c_p \right) \right] \\
\displaystyle\frac{\partial c_p}{\partial t} = D_p \left[\nabla \cdot \left(\nabla c_p + \frac{z_n\,e}{k_B\,T}\,c_p\,\nabla \phi \right) + \tilde{S_\sigma} \nabla \cdot \left(g_{np}\,c_p\,\nabla c_n + g_{pp}\,c_p\,\nabla c_p \right) \right] \\
-\nabla \cdot (\eps\,\nabla \phi) = \rho_0 + z_n\,e\,c_n + z_p\,e\,c_p
\end{cases}
\end{equation*}

Let $\eps=1$ and {$\rho_0=0$}  {in the above system}, we end up with the following two-component \textit{drift-diffusion system}
\begin{equation}\label{eqn: cn cp phi system with steric effect}
\begin{cases}
u_t=\nabla\cdot (d_1\,\nabla u+\vartheta_1\,u\,\nabla\phi)+\nabla\cdot (g_{11}\,u\,\nabla u+g_{12}\,u\,\nabla v), \quad & x\in\Omega, \quad t>0,\\
v_t=\nabla\cdot (d_2\,\nabla v+\vartheta_2\,v\,\nabla\phi)+\nabla\cdot (g_{21}\,v\,\nabla u+g_{22}\,v\,\nabla v), \quad & x\in\Omega, \quad t>0,\\
-\Delta \phi=\gamma_1\,u+\gamma_2\,v, \quad x\in\Omega, \quad t>0,
\end{cases}
\end{equation}
where $u=u(x,t)$ and $v=v(x,t)$ are assumed to be {positive} functions;  $d_1>0$ and $d_2>0$ are diffusion rates. Throughout this paper, we assume that $\vartheta_1$, $g_{11}$, $g_{12}$, $g_{21}$, $g_{22}$ and $\gamma_1$ are positive constants; $\vartheta_2$ and $\gamma_2$ are negative constants.
In this paper, we are concerned with stationary solutions to \eqref{eqn: cn cp phi system with steric effect}, i.e. with time-independent solutions to the following elliptic system
\begin{equation}\label{eqn: cn cp phi system with steric effect steady state}
\begin{cases}
0=\nabla\cdot (d_1\,\nabla u+\vartheta_1\,u\,\nabla\phi)+\nabla\cdot (g_{11}\,u\,\nabla u+g_{12}\,u\,\nabla v), \quad x\in\Omega,\\
0=\nabla\cdot (d_2\,\nabla v+\vartheta_2\,v\,\nabla\phi)+\nabla\cdot (g_{21}\,v\,\nabla u+g_{22}\,v\,\nabla v), \quad x\in\Omega,\\
-\Delta \phi=\gamma_1\,u+\gamma_2\,v, \quad x\in\Omega.
\end{cases}
\end{equation}
Using the fact that $\nabla (\log u)=\nabla u/u$, the first and second equations in \eqref{eqn: cn cp phi system with steric effect steady state} can be rewritten as 
\begin{equation}\label{eqn: cn cp phi system with steric effect steady state without phi}
\begin{cases}
0=\nabla\cdot \Big( u\, \nabla (d_1\,\log u+\vartheta_1\,\phi+g_{11}\,u+g_{12}\,v )\Big), \quad x\in\Omega,\\
0=\nabla\cdot \Big( v\, \nabla (d_2\,\log v+\vartheta_2\,\phi+g_{21}\,u+g_{22}\,v )\Big), \quad x\in\Omega.
\end{cases}
\end{equation}
It is readily seen that if we can find $u$, $v$ and $\phi$ satisfying the \textit{algebraic equations}
\begin{equation}\label{eqn: drift-diffusion DAEs algebraic eqns}
\begin{cases}
d_1\,\log u(x)+\vartheta_1\,\phi(x)+g_{11}\,u(x)+g_{12}\,v(x)=c_1, \quad x\in\Omega,\\
d_2\,\log v(x)+\vartheta_2\,\phi(x)+g_{21}\,u(x)+g_{22}\,v(x)=c_2, \quad x\in\Omega,
\end{cases}
\end{equation}
where $c_1$ and $c_2$ are constants, then such $u$, $v$ and $\phi$ automatically form a solution of \eqref{eqn: cn cp phi system with steric effect steady state without phi}.  A natural question arises as to whether \textit{any} solution of \eqref{eqn: cn cp phi system with steric effect steady state without phi} also satisfies \eqref{eqn: drift-diffusion DAEs algebraic eqns}.  It will be shown in \thref{prop: Equivalence of algebraic equations and differential equations} that the answer is indeed affirmative when certain appropriate boundary conditions are imposed on the solutions, i.e.
\begin{equation}\label{eqn: BC for F1}
u\,F_1\,\frac{\partial F_1}{\partial \nu}\leq0, \quad \text{a.e. on} \quad \partial\Omega,
\end{equation}
and
\begin{equation}\label{eqn: BC for F2}
v\,F_2\,\frac{\partial F_2}{\partial \nu}\leq0, \quad \text{a.e. on} \quad \partial\Omega,
\end{equation}
where $F_1:=d_1\,\log u+\vartheta_1\,\phi+g_{11}\,u+g_{12}\,v$ and $F_2:=d_2\,\log v+\vartheta_2\,\phi+g_{21}\,u+g_{22}\,v$. It is worth noticing that (\ref{eqn: BC for F1}) and (\ref{eqn: BC for F2}) are guaranteed when, for instance, $\frac{\partial F_i}{\partial \nu}=0 \ \ (i=1,2)$ on $\partial\Omega$, or the homogeneous Neumann boundary conditions hold:
\begin{equation}\label{eqn: zero Neumann BC}
\frac{\partial u}{\partial \nu}=\frac{\partial v}{\partial \nu}=\frac{\partial \phi}{\partial \nu}=0 \quad \text{on} \quad \partial\Omega.
\end{equation}
As a consequence, our problem now turns to establishing the existence and analyzing behavior of solutions to the \textit{differential algebraic equations} (DAEs): 
\begin{equation*}
\begin{cases}
d_1\,\log u+\vartheta_1\,\phi+g_{11}\,u+g_{12}\,v=c_1, \quad x\in\Omega,\\
d_2\,\log v+\vartheta_2\,\phi+g_{21}\,u+g_{22}\,v=c_2, \quad x\in\Omega,\\
-\Delta \phi=\gamma_1\,u+\gamma_2\,v, \quad x\in\Omega.
\end{cases}
\end{equation*}
The DAE approach is quite simple but efficient, and allows us to get a complete understanding of behavior of solutions to the steady-state PNP-steric equations of two ion species. The reader is referred to \cite{gavish2017poisson} for an analyze of PNP-steric equations via PNP-Cahn-Hilliard model. 

{The existence and basic properties of solutions $u=u(\phi)$ and $v=v(\phi)$ to the system \eqref{eqn: drift-diffusion DAEs algebraic eqns} for any parameters $d_1,d_2,\vartheta_1,g_{11},g_{12},g_{21},g_{22}, \gamma_1 >0$; $\vartheta_2, \gamma_2<0$; and $c_1,c_2 \in \mathbb{R}$ is verified in \thref{thm: Existence of solutions to AEs} and \thref{prop: alge eqns}. 

Moreover, under the hypothesis
\newline
\textbf{(H1)} $g_{11}\,g_{22}-g_{12}\,g_{21}\ge 0$, 
\newline
we will show in \thref{thm: Uniqueness of solutions to AEs} that \eqref{eqn: drift-diffusion DAEs algebraic eqns} admits unique solution $(u,v,\phi)$ in which $u$ and $v$ can be uniquely represented w.r.t. $\phi$, i.e. $u=u(\phi)$ and $v=v(\phi)$, and the curve $u(\phi),v(\phi)$ are of class $C^1$.} Note that \eqref{eqn: drift-diffusion DAEs algebraic eqns} is a system of nonlinear algebraic equations for which \textit{explicit solutions} expressed by the form $u=u(\phi)$ and $v=v(\phi)$ in general cannot be found. Due to $\textbf{(H1)}$ however, the solution {$u=u(\phi)$ and $v=v(\phi)$} of \eqref{eqn: drift-diffusion DAEs algebraic eqns} \textit{in implicit form} can be given \textit{uniquely}. With the aid of \thref{thm: Uniqueness of solutions to AEs}, we {arrive at} the following semilinear Poisson equation
\begin{equation}\label{eqn: semilinear Poisson eqn}
-\Delta \phi+G(\phi)=0,\quad x\in\Omega,
\end{equation}
where $G(\phi):=-\gamma_1\,u(\phi)-\gamma_2\,v(\phi)$. To establish existence of solutions of \eqref{eqn: semilinear Poisson eqn} under the zero Neumann boundary condition 
\begin{equation*}
\frac{\partial \phi}{\partial \nu}=0 \quad \text{on} \quad \partial\Omega,
\end{equation*}
{more delicate} properties of the nonlinearity $G=G(\phi)$ and the solution $u=u(\phi)$, $v=v(\phi)$ under \textbf{(H1)} are explored in \thref{prop: alge eqns H1}. 


\begin{theorem}[\textbf{Existence and uniqueness of $C^1$ solutions to the steady-state PNP-steric equations under $\textbf{(H1)}$}]\thlabel{thm: main H1}
Let $\Omega = (-1,1) \subset \mathbb{R}$. Assume that $\textbf{(H1)}$ holds. Then \eqref{eqn: cn cp phi system with steric effect steady state} coupled with the homogeneous Neumann boundary conditions \eqref{eqn: zero Neumann BC} has a unique solution $(u,v,\phi)=(u(\phi(x)),v(\phi(x)),\phi(x))$ for all $x \in \Omega$. Moreover, $u(\phi(x)),v(\phi(x))$ and $ \phi(x)$ are of class $C^2$ for all $x \in \bar{\Omega}$.


\end{theorem}
The {idea behind} the DAEs approach we use to obtain \thref{thm: main H1} is elementary. However, the result is remarkable in that only $\textbf{(H1)}$ is needed to ensure the uniqueness and $C^2$-smoothness of solutions to the elliptic system \eqref{eqn: cn cp phi system with steric effect steady state} under the homogeneous Neumann boundary conditions \eqref{eqn: zero Neumann BC}. 

On the other hand, if \textbf{(H1)} is violated, \eqref{eqn: drift-diffusion DAEs algebraic eqns} may admit multiple solutions. In Section~\ref{sec: Triple solutions}, we shall give a simple criterion to check whether the system  \eqref{eqn: drift-diffusion DAEs algebraic eqns} admits unique monotone $C^2$ solution; or unique monotone piecewise $C^2$ solution with a vertical tangent; or triple piecewise $C^2$ solution when \textbf{(H1)} fails. We would like to mention that the case of triple solution is very similar to the S-shaped solutions in \cite{ding2017computational,steinruck1989bifurcation}. More detailed analysis regarding the bifurcation of the elliptic system  \eqref{eqn: cn cp phi system with steric effect steady state} will be studied in the following-up work.

\section{Two species equations}\label{sec: Two species equations}


To begin with, we show that under the boundary conditions \eqref{eqn: BC for F1} and \eqref{eqn: BC for F2}, every solution of \eqref{eqn: cn cp phi system with steric effect steady state without phi} also solves \eqref{eqn: drift-diffusion DAEs algebraic eqns}, as mentioned in Section~\ref{Sec:model}.
\begin{proposition}[\textbf{Equivalence of algebraic and differential equations}] \thlabel{prop: Equivalence of algebraic equations and differential equations} {Under the boundary conditions  \eqref{eqn: BC for F1} and \eqref{eqn: BC for F2}, for any pair $(c_1,c_2)\in \mathbb{R}^2$, the systems of steady-state PNP-steric equations \eqref{eqn: cn cp phi system with steric effect steady state without phi}  is equivalent to the system of DAEs \eqref{eqn: drift-diffusion DAEs algebraic eqns}.}
\end{proposition}
\begin{proof} The fact that any solution to  \eqref{eqn: drift-diffusion DAEs algebraic eqns} solves \eqref{eqn: cn cp phi system with steric effect steady state without phi} is trivial. Now, we shall check that every solution $(u,v,\phi)$ to \eqref{eqn: cn cp phi system with steric effect steady state without phi} satisfies \eqref{eqn: drift-diffusion DAEs algebraic eqns}. Indeed, let $F$ as in Section~\ref{Sec:model}, i.e. $F_1:=d_1\,\log u+\vartheta_1\,\phi+g_{11}\,u+g_{12}\,v$, {the first equation {in} \eqref{eqn: cn cp phi system with steric effect steady state without phi} or} $\nabla\cdot(u\,\nabla F_1)$=0 gives $0=F_1\,\nabla\cdot(u\,\nabla F_1)=F_1\,\nabla u\cdot\nabla F_1+u\,F_1\Delta F_1.$
{On the other hand}, taking into account \eqref{eqn: BC for F1} and applying integration by parts yield
\begin{align*}
\int_{\Omega} u\,F_1\Delta F_1 \,dx
&=    -\int_{\Omega} \nabla(u\,F_1)\cdot \nabla F_1\,dx
        +\int_{\partial\Omega}u\,F_1\,\frac{\partial F_1}{\partial \nu}\,ds\notag\\
&\leq-\int_{\Omega} u\,|\nabla F_1|^2\,dx
        -\int_{\Omega} F_1\,\nabla u\cdot\nabla F_1\,dx. 
\end{align*}
Thus, we can confirm that $\int_{\Omega} u\,|\nabla F_1|^2\,dx\leq0$. Since $u(x)>0$ on $\Omega$, $F_1$ must be a constant independent of $x$ a.e. $\Omega$.

In a similar manner, we can prove that $F_2:=d_2\,\log v+\vartheta_2\,\phi+g_{21}\,u+g_{22}\,v$ is a constant independent of $x$  a.e. $\Omega$ (starting from $\nabla\cdot(v\,\nabla F_2)=0$). 

This completes the proof of \thref{prop: Equivalence of algebraic equations and differential equations}.
\end{proof}



\begin{remark}
When the zero Neumann boundary conditions $\displaystyle\frac{\partial u}{\partial \nu}=\frac{\partial v}{\partial \nu}=\frac{\partial \phi}{\partial \nu}=0$ on $\partial\Omega$ {are} considered, it is easy to see that these boundary conditions lead to $\displaystyle\frac{\partial F_1}{\partial \nu}=0$ and $\displaystyle\frac{\partial F_2}{\partial \nu}=0$ on $\partial\Omega$.  Thus, \eqref{eqn: BC for F1} and \eqref{eqn: BC for F2} hold true.


\end{remark}

We are now in the position to investigate existence of solutions to \eqref{eqn: drift-diffusion DAEs algebraic eqns}. 
\begin{theorem}[\textbf{Existence of solutions to \eqref{eqn: drift-diffusion DAEs algebraic eqns}}]\thlabel{thm: Existence of solutions to AEs}
For any $\phi \in \mathbb{R}$ {and any pair $(c_1,c_2)\in \mathbb{R}^2$}, there exists a solution $(u,v,\phi)$ to \eqref{eqn: drift-diffusion DAEs algebraic eqns}.
\end{theorem}

\begin{figure}[ht]
\center
\scalebox{0.55}[0.55]{\includegraphics{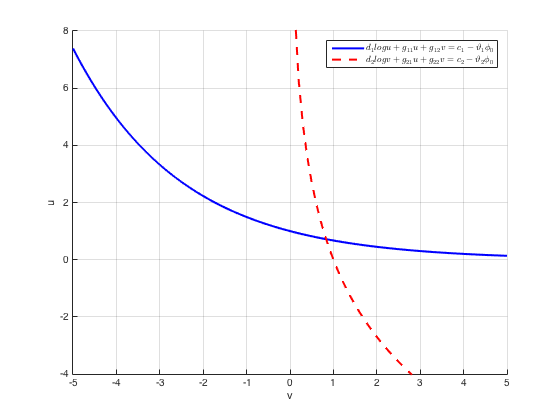}}
\caption{Example of the two curves satisfying (P1)-(P6) must have at least one intersection in the first quadrant.}
\label{fig: 2-species}
\end{figure}

\begin{proof}
Let $\phi_0 \in\mathbb{R}$ be fixed. It is sufficient to check that for {and any pair $(c_1,c_2)\in \mathbb{R}^2$} and for fixed $\phi_0$, there exists a solution $(u_0,v_0) \in \mathbb{R}^2$ to \eqref{eqn: drift-diffusion DAEs algebraic eqns}. In other words, our goal is to prove that the two curves ${(\mathcal{C}_1)}:  d_1\,\log u+g_{11}\,u+g_{12}\,v=c_1-\vartheta_1\,\phi_0$ and ${(\mathcal{C}_2)}: d_2\,\log v+g_{21}\,u+g_{22}\,v=c_2-\vartheta_2\,\phi_0$  {have} at least one intersection point in the first quadrant of the $v$-$u$ plane. 

To see this, we first rewrite $(\mathcal{C}_1)$ as follows
\[
v(u) = \frac{1}{g_{12}}\left(-d_1 \log u - g_{11}u + c_1 - \vartheta_1\phi_0\right).
\]
In the above equation, $v(u)$ can be understood as a smooth function w.r.t. $u$ for any $u>0$. Thus, we can differentiate it w.r.t. $u$ to obtain
\[
v'(u) = - \frac{1}{g_{12}}\left(\frac{d_1}{u} + g_{11}\right) <0.
\]
Since $v'(u)$ exists and negative for any $u>0$, its inverse function is also differentiable and negative. Therefore, differentiating $(\mathcal{C}_1):  d_1 \log u(v) +g_{11} u(v) + g_{12}v = c_1 - \vartheta_1 \phi_0$ w.r.t. $v$, we have
\[
u'(v) = -\frac{g_{12}\,u(v)}{d_1+g_{11}\,u(v)}<0.
\]
It is readily to see that the {graph of} $u=u(v)$ {satisfying} $d_1\,\log u+g_{11}\,u+g_{12}\,v=c_1-\vartheta_1\,\phi_0$ on the $v$-$u$ plane has the following propert{ies}:
\begin{itemize}
  \item[(P1)] as $v\rightarrow -\infty$, $u\rightarrow\infty$;
  \item[(P2)] as $v\rightarrow\infty$, $u\rightarrow 0^{+}$;
  \item[(P3)] $u=u(v)$ is decreasing in $v\in\mathbb{R}$.
\end{itemize} 

In the same manner, we get from the curve $(\mathcal{C}_2)$ that $v(u)>0$ (since the curve $(\mathcal{C}_2)$ is well-defined for $v>0$),
\begin{equation*}
v'(u)=-\frac{g_{21}\,v(u)}{d_2+g_{22}\,v(u)}<0.
\end{equation*}
and that the {graph of}  $v=v(u)$ satisfying $d_2\,\log v(u)+g_{21}\,u+g_{22}\,v(u)=c_2-\vartheta_2\,\phi_0$  on the $v$-$u$ plane enjoys the following {properties}:
\begin{itemize}
  \item[(P4)] as $u\rightarrow -\infty$, $v\rightarrow\infty$;
  \item[(P5)] as $u\rightarrow\infty$, $v\rightarrow 0^{+}$;
  \item[(P6)] $v=v(u)$ is decreasing in $u\in\mathbb{R}$.
\end{itemize} 

As a consequence, it follows from the {properties} of the {graphs} of the two curves {($\mathcal{C}_1$) and ($\mathcal{C}_2$)} that {these} two curves in the first quadrant of the $v$-$u$ plane intersect at least once (cf. Figure~\ref{fig: 2-species}). That is, given any $\phi_0\in\mathbb{R}$, we can find at least one solution $(u_0,v_0)\in \mathbb{R}^2$ which satisfies \eqref{eqn: drift-diffusion DAEs algebraic eqns}.
\end{proof}

Some basic properties of solutions $(u(\phi),v(\phi))$ to \eqref{eqn: drift-diffusion DAEs algebraic eqns} will be given in the following proposition.

\begin{proposition}[\textbf{Properties of solutions to \eqref{eqn: drift-diffusion DAEs algebraic eqns}}]\thlabel{prop: alge eqns}
Let $(u,v)=(u(\phi),v(\phi))$  be a pair of solutions to \eqref{eqn: drift-diffusion DAEs algebraic eqns}. Then the pair $(u,v)$ enjoys the following properties:
\begin{itemize}
 \item [(i)] \textbf{(Asymptotic behavior of $u(\phi)$ and $v(\phi)$)}
 \begin{itemize}
 \item[$\triangleright$] As $\phi \to \infty: u(\phi) \to 0, v(\phi) \to \infty$,
 \item[$\triangleright$] As $\phi \to -\infty: u(\phi) \to \infty, v(\phi) \to 0$.
 \end{itemize}

  \item [(ii)] \textbf{(Local uniqueness of $u(\phi)$ and $v(\phi)$)} For any fixed $(c_1,c_2) \in \mathbb{R}^2$, if $(u_1,{v}_1)$ and $({u}_2,{v}_2)$ both satisfy \eqref{eqn: drift-diffusion DAEs algebraic eqns} corresponding to $\phi_1$ and $\phi_2$, then for any $u_1 \le u \le u_2$, there exist unique  $v_2 \le v \le v_1$ and $\phi \in \mathbb{R}$ such that $(u,v,\phi)$ satisfies \eqref{eqn: drift-diffusion DAEs algebraic eqns}. 
  
  \item[(iii)] \textbf{(Formulas of $u'(\phi)$ and $v'(\phi)$)} If $u'(\phi)$ and $v'(\phi)$ exist, they can be defined uniquely by the following formulas
  \begin{equation}\label{eqn: u'(phi)}
  u'(\phi_0)=-\frac{\vartheta _1
   \left(\frac{d_2}{v(\phi_0
   )}+g_{22}\right)-g_{12}\, \vartheta _2}{\left(\frac{d_1}{u(\phi_0
   )}+g_{11}\right)
   \left(\frac{d_2}{v(\phi_0
   )}+g_{22}\right)-g_{12}\,
   g_{21}},
  \end{equation}
    \begin{equation}\label{eqn: v'(phi)}
  v'(\phi_0)=\frac{g_{21}\, \vartheta _1-\vartheta _2
   \left(\frac{d_1}{u(\phi_0
   )}+g_{11}\right)}{\left(\frac{d_1}{u(\phi_0
   )}+g_{11}\right)
   \left(\frac{d_2}{v(\phi_0
   )}+g_{22}\right)-g_{12}\,
   g_{21}}.
  \end{equation}

\end{itemize}

\end{proposition}
\begin{proof}
%
(i) follows \eqref{eqn: drift-diffusion DAEs algebraic eqns} immediately. Indeed, letting $\phi \to \infty$ in \eqref{eqn: drift-diffusion DAEs algebraic eqns}, we have
\[
d_1 \log u(\phi) + g_{11} u (\phi)+ g_{12} v(\phi) \to -\infty, \quad  d_2 \log v(\phi) +g_{21}u(\phi) +g_{22} v(\phi) \to \infty
\]
Since $u(\phi),v(\phi)>0$ for all $\phi$, we must have $\log u(\phi) \to -\infty$ as $\phi \to \infty$. This means that $u(\phi) \to 0$ as $\phi \to \infty$. Thus, $v(\phi) \to \infty$ as $\phi \to \infty$. Employing similar argument for the case $\phi \to -\infty$, we get (i).

For (ii), we consider the following two functions
\[
Q_1(x):= \vartheta_1d_2 \log x + (\vartheta_1g_{22} -\vartheta_2 g_{12}) x,
\]
and
\[
Q_2(y):= \vartheta_2 d_1 \log y + (\vartheta_2g_{11}-\vartheta_1 g_{21}) y + \vartheta_1 c_2 - \vartheta_2 c_1.
\]
It is obvious to see that $Q_1$ is monotone increasing, whilst $Q_2$ is  monotone decreasing. Moreover, $Q_1(v_1) = Q_2(u_1)$ and $Q_1(v_2)=Q_2(u_2)$. Thus, $Q_1([v_2,v_1])$ and $Q_2([u_1,u_2])$ have the same range. As a consequence, for any $u_1 \le u \le u_2$, there exists unique $v_2 \le v \le v_1$ such that $Q_1(v)=Q_2(u)$. 

To prove $(iii)$, we first differentiate the two equations in \eqref{eqn: drift-diffusion DAEs algebraic eqns} one by one with respect to $\phi$, and obtain two equations in which the unknowns can be viewed as $u'(\phi)$ and $v'(\phi)$. Solving them gives $u'(\phi)$ and $v'(\phi)$ as stated in $(iii)$.

\end{proof} 

%
%
%

\section{Unique solution under \textbf{(H1)}}\label{sec:H1}

The uniqueness of solutions to \eqref{eqn: drift-diffusion DAEs algebraic eqns}  under $\textbf{(H1)}$ will be {established} in the following theorem.
\begin{theorem}[\textbf{Uniqueness of  solutions to \eqref{eqn: drift-diffusion DAEs algebraic eqns} under $\textbf{(H1)}$}]\thlabel{thm: Uniqueness of solutions to AEs}
Assume $\textbf{(H1)}$. Then for any $\phi \in \mathbb{R}$, and any pair $(c_1,c_2) \in \mathbb{R}^2$, there exists a unique solution $(u,v,\phi) $ to \eqref{eqn: drift-diffusion DAEs algebraic eqns}, which can be represented implicitly as $(u,v)=(u(\phi),v(\phi))$ and $u(\phi), v(\phi)$ are $C^1$ functions.
\end{theorem}

\begin{figure}[ht]
\center
\scalebox{0.55}[0.55]{\includegraphics{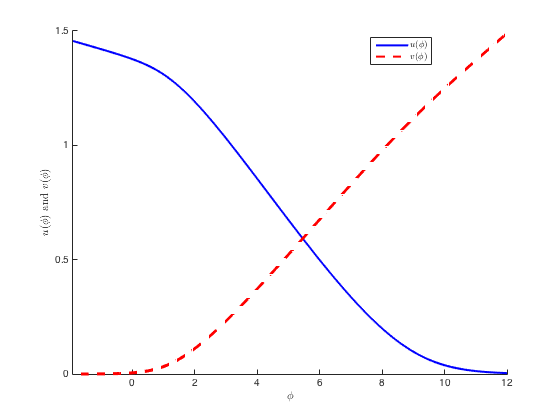}}
\caption{Unique solution under {\bf (H1)} (corresponds to parameters $d_1 = 1.0; 
d_2 = 2.0;
\vartheta_1 = 0.3;
\vartheta_2 = -4.0;
g_{11} = 7.0;
g_{12}= 8.0;
g_{21} = 9.0;
g_{22} = 33;
c_1 = 10;
c_2 = 2$).}
\label{fig:H1}
\end{figure}

\begin{proof}
The existence of solutions to \eqref{eqn: drift-diffusion DAEs algebraic eqns}  has been established in \thref{thm: Existence of solutions to AEs}. We now eliminate the possibility of non-uniqueness of solutions $(u,v)$ to \eqref{eqn: drift-diffusion DAEs algebraic eqns}  for a given $\phi_0 \in\mathbb{R}$ by contradiction. Suppose that, contrary to our claim, there exist in the first quadrant of the $v$-$u$ plane two distinct solutions $(u_1,v_1)$ and $(u_2,v_2)$ which satisfy \eqref{eqn: drift-diffusion DAEs algebraic eqns}  for a given $\phi_0\in\mathbb{R}$. {In the $v$-$u$ plane}, we consider the following functions 
\[
M_1: \R_+^2 \to \R: (u,v) \mapsto -\frac{g_{12}u}{d_1+g_{11}u}, \quad
M_2: \R_+^2 \to \R: (u,v) \mapsto -\frac{d_2+g_{22}v}{g_{21}v}
\]
and $V:\R_+^2 \to \R: (u,v) \mapsto M_1(u,v)-M_2(u,v)$.
It is worth noticing that $M_1(u,v)$ and $M_2(u,v)$ can be understood as the slope of {the curves $(\mathcal{C}_1): d_1 \log u +g_{11}u +g_{12}v =c_1-\vartheta_1 \phi_0 $ and $(\mathcal{C}_2): d_2 \log v +g_{21}u +g_{22}v =c_2-\vartheta_2 \phi_0$} at $(u,v)$, respectively. 
We consider the following three cases of {the quantity ${S}:=V(u_1,v_1)\cdot V(u_2,v_2)$.}
\begin{figure}[ht]
\center
\scalebox{0.55}[0.55]{\includegraphics{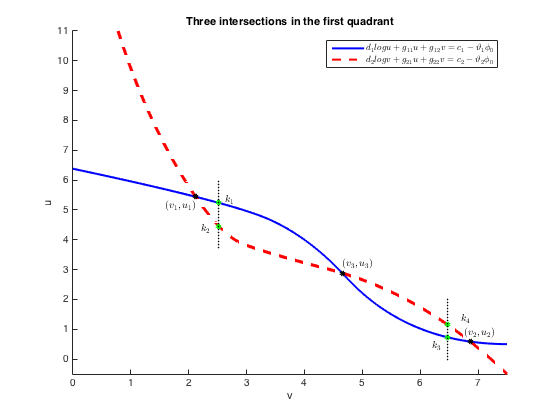}}
\caption{Example of the two curves satisfying (P1)-(P6) and obeying (S3).}
\label{fig:H1uv3}
\end{figure}
\begin{itemize}
\item[{(S1)}]  If $S=0$, then either $M_1(u_1,v_1)-M_2(u_1,v_1) =0$ or $M_1(u_2,v_2)-M_2(u_2,v_2)=0$. Without loss of generality, we assume that $M_1(u_1,v_1)-M_2(u_1,v_1) =0$. Taking into account definitions of $M_1$ and $M_2$, the fact that $M_1(u_1,v_1)=M_2(u_1,v_1)$ leads to 
\begin{equation*}
\frac{g_{12}\,u_1}{d_1+g_{11}\,u_1}=\frac{d_2+g_{22}\,v_1}{g_{21}\,v_1}.
\end{equation*}
It turns out that the last equation is equivalent to  
\begin{equation*}\label{eqn: eqn contradicting to (H1)}
\left(\frac{d_1}{u_1}+g_{11}\right)
\left(\frac{d_2}{v_1}+g_{22}\right)
=g_{12}\,g_{21},
\end{equation*}
which contradicts $\textbf{(H1)}$. 
\item[(S2)] If $S<0$, without loss of generality, we may assume that $M_1(u_1,v_1)-M_2(u_1,v_1)<0$ and $M_1(u_2,v_2)-M_2(u_2,v_2)>0$. Let {$h(t):= ((1-t)u_1+t\,u_2,(1-t)v_1+t\,v_2)$}, {then} the function ${V \circ h}: [0,1]\to \R$ is continuous and satisfies $V\circ h(0) <0$ and $V\circ h(1) >0$. By the Intermediate Value Theorem, there exists $t^{\ast} \in (0,1)$ for which $M_1(u^{\ast},v^{\ast})-M_2(u^{\ast},v^{\ast})=0$, where {$u^{\ast}:=(1-t^\ast)u_1 + t^\ast u_2$ and $v^\ast:=(1-t^\ast)v_1+t^\ast v_2$}.  Continuing {with} the argument in (S1) for $(u^{\ast},v^{\ast})$, we also get a contradiction to {\bf (H1)}.
\item[(S3)] If $S>0$, without loss of generality, we may assume that $M_1(u_1,v_1)<M_2(u_1,v_1)$ and $M_1(u_2,v_2)<M_2(u_2,v_2)$. For $\delta>0$ small enough, the vertical line at $v_1+\delta$ {must intersect} the curves $(\mathcal{C}_1)$ and $(\mathcal{C}_2)$ at $(v_1+\delta,k_1)$ and $(v_1+\delta,k_2)$ with $k_1>k_2$. Similarly, the vertical line at $v_2-\delta$ {must intersect} the curves $(\mathcal{C}_1)$ and $(\mathcal{C}_2)$ at $(v_2-\delta,k_3)$ and $(v_2-\delta,k_4)$ with $k_3<k_4$. Therefore, the two curves $(\mathcal{C}_1)$ and $(\mathcal{C}_2)$ {must intersect} once again at some point $(u_3,v_3)$ with $v_1<v_3<v_2$ and $u_1>u_3>u_2$ (see Figure~\ref{fig:H1uv3}). Moreover, it must hold that $M_1(u_3,v_3)>M_2(u_3,v_3)$. We then repeat the case (S2) to get a contradiction to {\bf (H1)}.
\end{itemize}
Thus, for given $\phi\in\mathbb{R}$, uniqueness of solutions to \eqref{eqn: drift-diffusion DAEs algebraic eqns} follows.  

The $C^1$-smoothness of $u(\phi)$ and $v(\phi)$ is guaranteed by the Implicit Function Theorem. Indeed, consider 
\[F: \mathbb{R}_+^2 \times \mathbb{R} \to \mathbb{R}^2: (x,y,z) \mapsto \begin{pmatrix} d_1 \log x + \vartheta_1 z + g_{11}x+g_{12}y-c_1 \\ d_2 \log y + \vartheta_2 z + g_{21}x + g_{22}y-c_2\end{pmatrix}.\]
It can be seen that $F$ is a continuously differentiable function and $F(u(\phi),v(\phi),\phi) = 0$ for $\phi \in \mathbb{R}$. Moreover, the Jacobian matrix of $F$
\[
J_{F}(x,y):= \begin{pmatrix}\displaystyle\frac{\partial F_1}{\partial x}  & \displaystyle\frac{\partial F_1}{\partial y} \\ \displaystyle\frac{\partial F_2}{\partial x} & \displaystyle\frac{\partial F_2}{\partial y}  \end{pmatrix} = \begin{pmatrix} \displaystyle\frac{d_1}{x} + g_{11} & g_{12} \\ g_{21} & \displaystyle\frac{d_2}{y} + g_{22}\end{pmatrix}
\]
is positive for all pair $(x,y)\in \mathbb{R}_+^2$ when \textbf{(H1)} holds true. This means that  the Jacobian matrix of $F$ is invertible at each point $(u(\phi),v(\phi),\phi)$. The Implicit Function Theorem then implies that at every $\phi_0 \in \mathbb{R}$, there exists  an open set $U \subset \mathbb{R}$ containing $\phi_0$, such that there exists a unique continuously differentiable function $g:U \to \mathbb{R}_+^2$ such that $g(\phi_0)= (u(\phi_0),v(\phi_0))$ and $F(g(\phi),\phi) = 0$ for all $\phi \in U$. Thus, we obtain $C^1$ smoothness of the solution $(u,v)=(u(\phi),v(\phi))$ for $\phi \in \mathbb{R}$. 

Numerical solution $(u(\phi),v(\phi))$ of \eqref{eqn: drift-diffusion DAEs algebraic eqns} when the coupling parameters satisfy \textbf{(H1)} is demonstrated in Figure~\ref{fig:H1}.
\end{proof}

Apart from \thref{prop: alge eqns}, more important properties of solutions to \eqref{eqn: drift-diffusion DAEs algebraic eqns} under \textbf{(H1)} are investigated in the next proposition.
\begin{proposition}[\textbf{Properties of solutions to \eqref{eqn: drift-diffusion DAEs algebraic eqns} under $\textbf{(H1)}$}]\thlabel{prop: alge eqns H1}
Assume that $\textbf{(H1)}$ holds.
For any $(c_1,c_2)\in \mathbb{R}^2$, and any $\phi \in \mathbb{R}$, the system \eqref{eqn: drift-diffusion DAEs algebraic eqns} is uniquely solvable by the implicit functions $(u,v)=(u(\phi),v(\phi))$, where $u(\phi), v(\phi)$ are of class $C^2$. Moreover, $u=u(\phi)$ is monotonically decreasing in $\phi\in\mathbb{R}$, while $v=v(\phi)$ is monotonically increasing in $\phi\in\mathbb{R}$. In addition,
\begin{equation}\label{eqn: upper bound for G'(phi)}
\gamma_1\,u'(\phi)+\gamma_2\,v'(\phi)\leq -\kappa,
\end{equation}
for some constant $\kappa>0$ independent of $\phi$.

\end{proposition}
\begin{proof}
The fact that $u(\phi),v(\phi)$ are of class $C^1$ is due to \thref{thm: Uniqueness of solutions to AEs}. Thanks to \textbf{(H1)}, we know that
\[
\left(\frac{d_1}{u}+g_{11}\right)\left(\frac{d_2}{v}+g_{22}\right)-g_{12}g_{21} >0
\]
for all pair of positive numbers $(u,v)$. In view of \eqref{eqn: u'(phi)} and \eqref{eqn: v'(phi)}, we immediately get $u'(\phi)<0$ and $v'(\phi)>0$ for $\phi\in\mathbb{R}$. The $C^2$-smoothness of $u(\phi)$ and $v(\phi)$ is then obtained by differentiating \eqref{eqn: u'(phi)} and \eqref{eqn: v'(phi)} w.r.t. $\phi$.

Fix $-\infty < \phi_1 < \phi_2 < \infty$. By the monotonicity of $u(\phi)$ and $v(\phi)$, we have
\begin{itemize}
\item If $\phi < \phi_1$, then $u(\phi) > u_1 :=u(\phi_1)$ and $v(\phi) < v_1:=v(\phi_1)$.
\item If $\phi > \phi_2$, then $u(\phi) < u_2:=u(\phi_2)$ and $v(\phi) > v_2:=v(\phi_2)$.
\item If $\phi_1 \le \phi \le \phi_2$, then $u_1 > u(\phi)> u_2$ and $v_1 <v(\phi)< v_2$.
\end{itemize}
Denote by $x(\phi):=\frac{d_1}{u(\phi)}+g_{11} >0; y(\phi):= \frac{d_2}{v(\phi)}+g_{22} >0; a:=\gamma_2 \vartheta_2 >0; b:= \gamma_1 \vartheta_1>0.$ Let 
$$\kappa:= \frac{1}{2} \min\left\{\frac{a}{\frac{d_2}{v_1} + g_{22}},\frac{a}{\frac{d_2}{v_2} + g_{22}}, \frac{b}{\frac{d_1}{u_1}+g_{11}}, \frac{b}{\frac{d_1}{u_2}+g_{11}}\right\}>0.$$
It is sufficient to check that $ax(\phi) + by (\phi) - \kappa x(\phi) y(\phi) >0$ for all $\phi \in \mathbb{R}$. Indeed,
\begin{itemize}
\item If $\phi < \phi_1$, the fact that $x(\phi) < \frac{d_1}{u_1}+g_{11}$ implies $\kappa x(\phi) < b$. Thus, 
$$ax(\phi) + y(\phi) (b-\kappa x(\phi)) >0.$$
\item If $\phi > \phi_2$, then $y(\phi) < \frac{d_2}{v_2} + g_{22}$. Hence, $\kappa y(\phi) < a$. We have
\[
x(\phi) (a - \kappa y(\phi)) + b y(\phi) >0.
\]
\item If $\phi_1 \le \phi \le \phi_2$, it holds that $x(\phi) \le \frac{d_1}{u_2}+g_{11}$ and $\kappa x(\phi) < b$. We arrive at
\[
ax(\phi) + y(\phi) (b-\kappa x(\phi)) >0.
\]
\end{itemize}
\end{proof}

\thref{thm: Uniqueness of solutions to AEs} and \thref{prop: alge eqns},  \thref{prop: alge eqns H1} inspire us to consider the Neumann problem for the semilinear Poisson equation \eqref{eqn: semilinear Poisson eqn}, i.e.
\begin{equation}\label{eqn: Dirichlet problem for the semilinear Poisson equation}
\begin{cases}
-\Delta \phi+G(\phi)&=0\quad  \text{in}\quad \Omega,\\
\hspace{8mm}\frac{\displaystyle\partial \phi}{\displaystyle\partial \nu}& =0 \quad   \text{on}\quad \partial\Omega.
\end{cases}
\end{equation}
Here $G(\phi):=-\gamma_1\,u(\phi)-\gamma_2\,v(\phi)$, and $(u(\phi),v(\phi))$ is the unique solution to \eqref{eqn: drift-diffusion DAEs algebraic eqns} defined in \thref{{thm: Uniqueness of solutions to AEs}}. Due to \thref{prop: alge eqns H1}, the nonlinearity $G: \mathbb{R} \to \mathbb{R}$ is of class $C^2$. Moreover, \eqref{eqn: upper bound for G'(phi)} guarantees a positive constant $\kappa$ independent of $t$ such that $G'(t) >\kappa >0$ for all $t\in \mathbb{R}$. It is worth noticing that this property implies that $G$ is strictly monotone increasing, i.e.
\[
\left(G(s)-G(t)\right)(s-t) > 0 \quad \mbox{for all} \quad s,t\in \mathbb{R} \quad \mbox{with} \quad s\ne t.
\]
%
%

In the following theorem, the existence and uniqueness of a weak solution $\phi \in H^1(\Omega)$ to \eqref{eqn: Dirichlet problem for the semilinear Poisson equation} is considered.

\begin{theorem}[\textbf{Existence and uniqueness of solutions to Poission equations with $C^1$  nonlinearity}]\thlabel{thm: Existence of solns to Poission eqns with cont nonlinearity}
Let $\Omega \subset \mathbb{R}^d, (d \ge 1)$ be a bounded domain with smooth boundary. Assume that the nonlinearity $G: \mathbb{R} \to \mathbb{R}$ satisfies the following properties:
\begin{itemize}
  \item [$\textbf{(G1)}$] $t\mapsto G(t)$ is of class $C^1$ for all $t \in \mathbb{R}$.
\item [$\textbf{(G2)}$] there exist a constant $\kappa >0$ independent of $t$ such that $G'(t) > \kappa >0$ for all $t\in \mathbb{R}$.
\end{itemize}
Then the semilinear Poisson equation under homogeneous Neumann boundary condition \eqref{eqn: Dirichlet problem for the semilinear Poisson equation} admits a unique solution $\phi \in H^1(\Omega)$ such that $G(\phi) \in L^1(\Omega)$, $G(\phi)\phi \in L^1(\Omega)$ and
\begin{equation}\label{eqn:variational-Poisson-cont}
\int_{\Omega} \nabla \phi  \cdot \nabla v + \int_{\Omega} G(\phi)v = 0,
\end{equation}
for all $v\in H^1(\Omega)$. 
\end{theorem}
\begin{proof} We shall modify the proof in \cite{looker2006semilinear,webb1980boundary} to get the desired result. First, we consider the following regularized variational equation for all $\eps >0$ and all $n \in \mathbb{N}$
\begin{equation}\label{eqn:regularized-variational-Poisson-cont}
\eps \int_{\Omega} uv + \int_{\Omega} \nabla u \cdot \nabla v + \int_{\Omega} G_n(u)v = 0, \quad \mbox{ for all } \quad u,v \in H^1(\Omega).
\end{equation}
Here, the truncation of the nonlinearity was introduced in \cite{looker2006semilinear,webb1980boundary} to control the growth of the nonlinearity
\begin{equation}\label{def:Gn}
G_n(t) := \begin{cases} G(t),  &\quad \mbox{if} \quad |G(t)| \le n, \\ n\displaystyle\frac{G(t)}{|G(t)|},  &\quad \mbox{otherwise}.\end{cases}
\end{equation}
It can be check that for all $n \in \mathbb{N}$, $G_n$ is monotone increasing, i.e.
\[
(G_n(s)-G_n(t))(s-t) \ge 0 \quad \mbox{for all} \quad s,t \in \mathbb{R}.
\]

Moreover, \thref{le:Sn-pseudo} also shows that the operator $S_n: H^1(\Omega) \to (H^1(\Omega))^*$ satisfying
\[
\left<S_nu,v\right>  = \int_{\Omega} G_n(u)v, \quad u,v \in H^1(\Omega),
\]
 is well-defined for all $n\in \mathbb{N}$ and is pseudomonotone for $n > |G(0)|$.
 
 On the other hand, it follows \cite{looker2006semilinear} that the operator $L_\eps: H^1(\Omega) \to (H^1(\Omega))^*$ defined by
 \[
\left<L_\eps u,v\right> := \eps\int_\Omega uv + \int_\Omega \nabla u \cdot \nabla v, \quad u,v \in H^1(\Omega),
\]
 is linear, bounded, coercive, monotone increasing, symmetric and strongly continuous. Thus, the operator $L_\eps + S_n$ is bounded, pseudomonotone and coercive. Indeed, thanks to \thref{le:Sn-pseudo}, for $n>|G(0)|$,
 \[
\left<(L_\eps + S_n)u,u\right> \ge \left<L_\eps u,u\right> + \int_\Omega G(0)u \ge \left<L_\eps u,u\right> - |G(0)| |\Omega|^{\frac{1}{2}}\|u\|_{H^1(\Omega)}.
\]
Hence,
\[
\lim_{\|u\|_{H^1(\Omega)} \to \infty} \frac{\left<(L_\eps+S_n)u,u\right>}{\|u\|_{H^1(\Omega)}} = +\infty
\]
by the coercivity of $L_\eps$. 

With the help of \thref{thm:existence-pseudo}, we can conclude that for all $\eps>0$ and all $n>|G(0)|$, the system \eqref{eqn:regularized-variational-Poisson-cont} admits a solution $u^\eps_n \in H^1(\Omega)$.

We then show that the sequence $u^\eps_n$ weakly converges to some element $u^\eps$ in $H^1(\Omega)$ as $n \to \infty$. Moreover, $u^\eps$ satisfies
\begin{equation}\label{eqn:regularized-variational-Poisson-cont-eps}
\eps \int_{\Omega} u^\eps v + \int_{\Omega} \nabla u^\eps \cdot \nabla v + \int_{\Omega} G(u^\eps)v = 0, \quad \mbox{ for all } \quad v \in H^1(\Omega).
\end{equation}
In fact, due to the coercivity of $L_\eps$, there exists a constant $c>0$ independent of $n$ such that
\begin{align*}
0 \le c\|u_n^\eps\|^2_{H^1(\Omega)} &\le \left< L_\eps u_n^\eps,u_n^\eps\right> = - \left<S_n u_n^\eps,u_n^\eps\right> = - \int_{\Omega} G_n(u_n^\eps)u_n^\eps  \le -\int_\Omega G(0)u_n^\eps\\
&\le  |G(0)|\displaystyle |\Omega|^{\frac{1}{2}} \|u^\eps_n\|_{L^2(\Omega)} \le |G(0)||\Omega|^{\frac{1}{2}} \|u^\eps_n\|_{H^1(\Omega)}.
\end{align*}
The above inequalities yields that the sequence $u_n^\eps$ is uniformly bounded in $H^1(\Omega)$ by a constant independent of $n$. The the fact that the Hilbert space $H^1(\Omega)$ is reflexive then implies (up to subsequence) $u_n^\eps \rightharpoonup u^\eps$ weakly  in $H^1(\Omega)$. Also, (up to subsequence) $u_n^\eps \to u^\eps$ strongly in $L^2(\Omega)$. We then arrive at
\begin{equation}\label{eq:ueps}
\eps \int_\Omega u^\eps_n v + \int_\Omega \nabla u^\eps_n \cdot \nabla v \to \eps \int_\Omega u^\eps v + \int_\Omega \nabla u^\eps \cdot \nabla v \quad \mbox{as} \quad n \to \infty.
\end{equation}
On the other hand, the fact that $u^\eps_n$ satisfies \eqref{eqn:regularized-variational-Poisson-cont} leads to
\begin{align*}
0\le - \int_\Omega G_n(u^\eps_n)u^\eps_n &= \eps \int_\Omega |u^\eps_n|^2 + \int_\Omega |\nabla u^\eps_n|^2 \le (\eps + 1) \|u^\eps_n\|^2_{H^1(\Omega)} \le C^\eps,
\end{align*}
for some constant $C^\eps$ depending on $\eps$ and independent of $n$. By \thref{le:Gn-and-G}, $G_n(u_n^\eps) \to G(u^\eps)$ strongly in $L^1(\Omega)$ as $n \to \infty$. This yields that
\begin{equation}\label{eq:ueps-2}
\int_\Omega G_n (u^\eps_n)v \to \int_\Omega G(u^\eps)v \quad \mbox{ as } \quad n \to \infty, \quad \mbox{ for all } \quad v \in H^1(\Omega).
\end{equation}
Combining \eqref{eq:ueps} and \eqref{eq:ueps-2}, we can conclude that $u^\eps$ satisfies \eqref{eqn:regularized-variational-Poisson-cont}.

We end the existence part by showing that $u^\eps$ weakly converges to some $u$ in $H^1(\Omega)$ satisfying \eqref{eqn:variational-Poisson-cont}. Indeed, taking into account \textbf{(G2)}, we arrive at
\[
0 \le \eps \|u^\eps\|^2_{L^2(\Omega)} + \|\nabla u^\eps\|^2_{L^2(\Omega)} = - \int_\Omega G(u^\eps)u^\eps  \le - \int_\Omega G(0)u^\eps - \kappa \int_\Omega |u^\eps|^2.
\]
Thus, 
\[
\eps \|u^\eps\|^2_{L^2(\Omega)} + \|\nabla u^\eps\|^2_{L^2(\Omega)} + \kappa \|u^\eps\|^2_{L^2(\Omega)} \le -\int_\Omega G(0)u^\eps \le |G(0)||\Omega|^{\frac{1}{2}} \|u^\eps\|_{L^2(\Omega)}.
\]
This gives the uniformly bounded (independent of $\eps$) of $u^\eps$ in $H^1(\Omega)$, that is
\[
\min\{\kappa,1\} \|u^\eps\|_{H^1(\Omega)} \le |G(0)||\Omega|^{\frac{1}{2}},
\]
which implies that (up to subsequence) $u^\eps$ weakly converges to some $u$ in $H^1(\Omega)$. Also, (up to subsequence) we can assume that $u^\eps$ strongly converges to $u$ in $L^2(\Omega)$ and pointwisely converges to $u$ a.e. $\Omega$. As a consequence, for all $v \in H^1(\Omega)$
\[
\eps \int_\Omega u^\eps v + \int_\Omega \nabla u^\eps \cdot \nabla v \to \int_\Omega \nabla u \cdot \nabla v \quad \mbox{ as } \quad \eps \to 0.
\]
On the other hand, for small $\eps$ (i.e. $0< \eps \le 1$)
\[
0 \le - \int_\Omega G(u^\eps)u^\eps = \eps \int_\Omega |u^\eps|^2 + \int_\Omega |\nabla u^\eps|^2  \le \|u^\eps\|^2_{H^1(\Omega)} \le C,
\]
for some constant $C>0$ independent of $\eps$. Applying \thref{le:Gn-and-G}, we have that $G(u^\eps)$ strongly converges to $G(u)$ in $L^1(\Omega)$. This implies that $u$ satisfies \eqref{eqn:variational-Poisson-cont}.

We now prove the uniqueness by contradiction. Assume that $\phi_1$ and $\phi_2$ both satisfy \eqref{eqn:variational-Poisson-cont}. The strict monotonicity of $G$ gives
\[
\int_\Omega |\nabla (\phi_1-\phi_2)|^2 + \int_\Omega (G(\phi_1)-G(\phi_2))(\phi_1-\phi_2) >0,
\]
for $\phi_1 \ne \phi_2$. Since both $\phi_1$ and $\phi_2$ satisfy \eqref{eqn:variational-Poisson-cont}, they must satisfy
\[
\int_\Omega \nabla (\phi_1-\phi_2)\cdot \nabla v + \int_\Omega (G(\phi_1)-G(\phi_2))v = 0,
\]
for all $v \in H^1(\Omega)$. Choosing $v=\phi_1-\phi_2$ in the above equality, we then get a contradiction.
\end{proof}

We are now in a position to prove  \thref{thm: main H1}.
\begin{proof}[Proof of \thref{thm: main H1}]
Thanks to \thref{thm: Existence of solutions to AEs} and \thref{thm: Uniqueness of solutions to AEs}, the algebraic system \eqref{eqn: drift-diffusion DAEs algebraic eqns} admits a unique solution $(u(\phi),v(\phi),\phi)$, where $u(\phi)$ and $v(\phi)$ are of class $C^1$. Denote by $G(\phi):=-\gamma_1 u(\phi) - \gamma_2 v(\phi)$, then $G$ satisfies conditions \textbf{(G1)} and \textbf{(G2)} due to \thref{prop: alge eqns H1}. Upon using \thref{thm: Existence of solns to Poission eqns with cont nonlinearity}, we establish the existence of a unique solution $\phi \in H^1(\Omega)$ to the semilinear Poisson equation under homogeneous Neumann boundary \eqref{eqn: Dirichlet problem for the semilinear Poisson equation}. 

For $\Omega = (-1,1) \subset \mathbb{R}$, we can apply Corollary 8.11 in  \cite{brezis2010functional} to see that $u \circ \phi$ and $v \circ \phi$ actually belong to $H^1(\Omega)$. As a consequence, $G(\phi(\cdot)) \in H^1(\Omega)$. Hence, from $ \phi''(x) = G(\phi(x))$ we have $\phi \in H^3(\Omega)$.

Notice that $H^1(\Omega)$ consists of equivalence class of functions agreeing a.e. on $\Omega$. For each $f \in H^1(\Omega)$, there exists a unique continuous representative that agrees with $f$ a.e. $\Omega$ (see Theorem 8.2 of \cite{brezis2010functional}). Thus, we can assume that $\phi''(x)$ is a continuous function on $\overline{\Omega}$. This yields that $\phi \in C^2(\overline{\Omega})$. Therefore, $u(\phi(x)), v(\phi(x))$ are also of class $C^2$ due to the chain rule. Besides, the chain rule also implies
\[
\frac{\partial u(\phi(x))}{\partial \nu} = u'(\phi(x))\frac{\partial \phi(x)}{\partial \nu} =0, \quad \frac{\partial v(\phi(x))}{\partial \nu} = v'(\phi(x))\frac{\partial \phi(x)}{\partial \nu} =0
\]
for all $x \in \partial \Omega$, which means $u(\phi(x))$ and $v(\phi(x))$ both satisfy the homogeneous Neumann boundary conditions for $u$ and $v$.

Since the homogeneous Neumann boundary conditions \eqref{eqn: zero Neumann BC} guarantees \eqref{eqn: BC for F1} and \eqref{eqn: BC for F2}, we can employ \thref{prop: Equivalence of algebraic equations and differential equations} to complete the proof of \thref{thm: main H1}.
\end{proof}

\section{Bifurcation when \textbf{(H1)} is violated}\label{sec: Triple solutions}

{Throughout this section we consider the system \eqref{eqn: drift-diffusion DAEs algebraic eqns} when \textbf{(H1)} does not hold, that is, when $g_{11}\,g_{22}-g_{12}\,g_{21}< 0$}. 

\begin{figure}[ht]
\center
\scalebox{0.55}[0.55]{\includegraphics{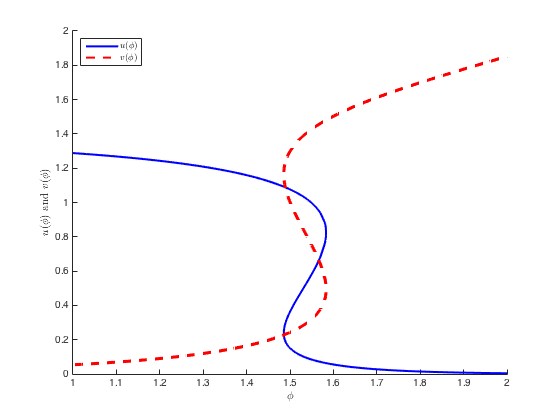}}
\caption{\textbf{{A} {t}riple solution when $\textbf{(H1)}$ is violated} (corresponds to parameters $d_1 = 1.0;
d_2 = 2.0;
\vartheta_1 = 0.3;
\vartheta_2 = -4.0;
g_{11} = 7.0;
g_{12} = 8.0;
g_{21} = 9.0;
g_{22} = 33/7;
c_1 = 10;
c_2 = 2)$. 
}
\label{fig: u=u(phi)}
\end{figure}

We can see that if \textbf{(H1)} fulfills, the quantity
\begin{equation*}
I(u,v):=\left(\frac{d_1}{u} +g_{11}\right)\left(\frac{d_2}{v}+g_{22}\right) - g_{12}g_{21} 
\end{equation*}
never vanishes for all pair $(u,v) \in \mathbb{R}_+^2$. This is the key point for \thref{thm: Uniqueness of solutions to AEs} to prove the uniqueness and $C^1$ smoothness of the solutions $(u(\phi),v(\phi))$ to the system \eqref{eqn: drift-diffusion DAEs algebraic eqns} for all $\phi \in \mathbb{R}$. Moreover, thanks to \textbf{(H1)}, $I(u(\phi),v(\phi)) \ne 0$ along the curves $u(\phi)$, $v(\phi)$,  and their slopes $u'(\phi)$, $v'(\phi)$ are assigned merely finite value along these curves. Therefore, it motivates us to investigate such points $(u(\phi),v(\phi))$ satisfying \eqref{eqn: drift-diffusion DAEs algebraic eqns} at which $|u'(\phi)|=|v'(\phi)|=\infty$ in order to have \textbf{(H1)} violated, i.e. points of the graphs of $u = u(\phi)$ and $v=v(\phi)$ having vertical tangent  as shown in Figure~\ref{fig: u=u(phi)}.

In the following, our aim is to find such points $(u,v)$ satisfying both \eqref{eqn: drift-diffusion DAEs algebraic eqns} and $I(u,v)=0$ by solving the {following algebraic} equations {for the unknowns $(u,v,\phi)$}:
\begin{equation}\label{eqn: algebraic eqns for triple solns}
\begin{cases}
d_1\,\log u+\vartheta_1\,\phi+g_{11}\,u+g_{12}\,v=c_1,\\
d_2\,\log v+\vartheta_2\,\phi+g_{21}\,u+g_{22}\,v=c_2,\\
\left(\displaystyle\frac{d_1}{u}+g_{11}\right)\left(\displaystyle\frac{d_2}{v}+g_{22}\right)-g_{12}\,g_{21}=0.
\end{cases}
\end{equation}
From the third equation in \eqref{eqn: algebraic eqns for triple solns}, we obtain 
\[
v\left(u(g_{12}g_{21}-g_{11}g_{22})-d_1g_{22}\right)=d_2(d_1+ug_{11})>0.
\]
Thus, if $u$ satisfies \eqref{eqn: algebraic eqns for triple solns}, then $u(g_{12}g_{21}-g_{11}g_{22})-d_1g_{22}>0$, and we can write
\begin{equation}\label{eqn: v=v(u)}
v=\displaystyle\frac{d_2\,\left(d_1+g_{11}\,u\right)}{u\,(g_{12}\,g_{21}-g_{11}\,g_{22})-d_1\, g_{22}}.
\end{equation}
Multiplying the first equation in \eqref{eqn: algebraic eqns for triple solns} by $\vartheta_2$ and the second equation in \eqref{eqn: algebraic eqns for triple solns} by $\vartheta_1$, we obtain two equations. Using \eqref{eqn: v=v(u)} and subtracting of one of the two equations from the other give $\sigma(u)=0$, where 
\begin{align*}
      \sigma(u)
:=&\frac{1}{\vartheta_1}\,\left(
          c_1-d_1\,\log u-g_{11}\,u
          -g_{12}\,\frac{d_2\,\left(d_1+g_{11}\,u\right)}
                  {u\,(g_{12}\,g_{21}-g_{11}\,g_{22})-d_1\, g_{22}}
     \right)-\notag\\
&\frac{1}{\vartheta_2}\,\bigg(
          c_2-d_2\,
          \log \left( \frac{d_2\,\left(d_1+g_{11}\,u\right)}
                  {u\,(g_{12}\,g_{21}-g_{11}\,g_{22})-d_1\, g_{22}}
                  \right)
\notag\\
&  \hspace{32mm}      
        -g_{21}\,u
        -g_{22}\,
          \frac{d_2\,\left(d_1+g_{11}\,u\right)}
                  {u\,(g_{12}\,g_{21}-g_{11}\,g_{22})-d_1\, g_{22}}
     \bigg).
\end{align*}

\begin{figure}[ht]
\center
\scalebox{0.55}[0.55]{\includegraphics{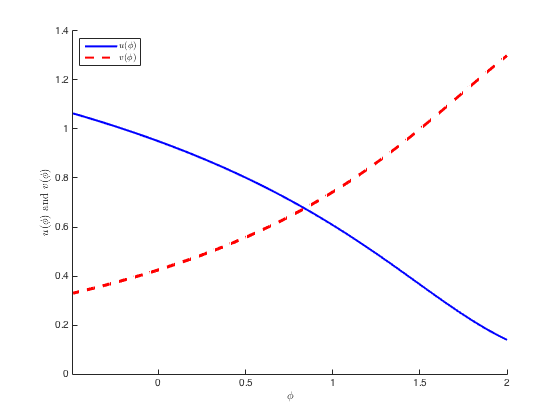}}
\caption{\textbf{Unique monotone solution when $\textbf{(H1)}$ is violated} (corresponds to parameters $d_1 = 1.0;
d_2 = 10.0;
\vartheta_1 = 0.3;
\vartheta_2 = -4.0;
g_{11} = 7.0;
g_{12} = 8.0;
g_{21} = 9.0;
g_{22} = 33/7;
c_1 = 10;
c_2 = 2)$. 
}
\label{fig:H22}
\end{figure}

\begin{figure}[ht]
\center
\scalebox{0.55}[0.55]{\includegraphics{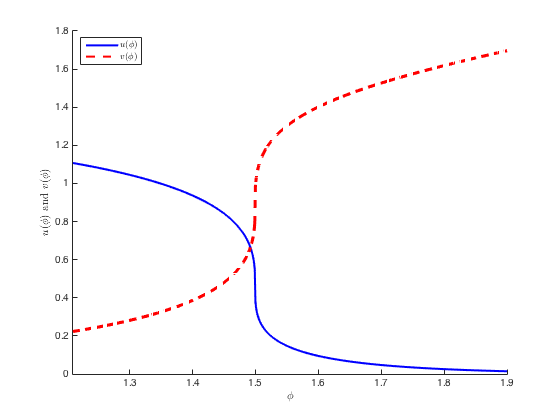}}
\caption{\textbf{Unique monotone solutions with vertical tangents when $\textbf{(H1)}$ is violated} (corresponds to parameters $d_1 = 1.0;
{d_2 = 2.784085121596521};
\vartheta_1 = 0.3;
\vartheta_2 = -4.0;
g_{11} = 7.0;
g_{12} = 8.0;
g_{21} = 9.0;
g_{22} = 33/7;
c_1 = 10;
c_2 = 2)$. 
}
\label{fig:H23}
\end{figure}

Now the question remains to determine the {graph} of {$\sigma=\sigma(u)$ on the $u$-$\sigma$ plane}. To this end, we observe that $\sigma=\sigma(u)$ {is defined} {for} $u>u^{\ast}$, where $u^{\ast}:=\displaystyle\frac{d_1\,g_{22}}{g_{12}\,g_{21}-g_{11}\,g_{22}}>0$. Also, it is readily verified that
\begin{equation}\label{eqn: asym behav of sigma(u) at infty}
\lim_{u\rightarrow (u^{\ast})^{+}} \sigma(u)=
\lim_{u\rightarrow \infty}             \sigma(u)=-\infty.
\end{equation}
To {determine} {the critical} points of $\sigma(u)=0$, we {find}
\begin{equation}\label{eqn: sigma'(u)}
\sigma'(u)=\frac{p(u)\,
\left(
u\,(g_{21}\,\vartheta_1-g_{11}\,\vartheta_2
\right)
-d_1\,\vartheta_2)}
{\vartheta _1\,\vartheta _2\,u\,
   \left(d_1+g_{11}\,u\right)\,
   \left(
   u\,(g_{12}\,g_{21}-g_{11}\,g_{22})-d_1\, g_{22}
   \right)^2},
\end{equation}
where $p(u):=k_3\,u^3+k_2\,u^2+k_1\,u+k_0$, and 
\begin{align}\nonumber
k_3&:= g_{11} \left(g_{12}\,g_{21}-g_{11}\, g_{22}\right)^2,\notag\\
k_2&:=d_1 \left(g_{12}\, g_{21}-3 \,g_{11}\,g_{22}\right) \left(g_{12}\,g_{21}-g_{11}\, g_{22}\right),\notag\\
k_1&:={{-}}d_1 \left(d_2\, g_{21}\, g_{12}^2+2\, d_1\,g_{21}\, g_{22}\, g_{12}-3 \,d_1\,g_{11}\, g_{22}^2\right), \notag\\
k_0&:=d_1^3\, g_{22}^2.\notag
\end{align}
We remark that the denominator of $\sigma'(u)$ {in \eqref{eqn: sigma'(u)}} {cannot be} $0$ since $u>u^{\ast}$. On the other hand, the numerator of {$\sigma'(u)$ in} \eqref{eqn: sigma'(u)} may admit up to four roots: 
\[
\displaystyle\frac{d_1\,\vartheta_2}{g_{21}\,\vartheta_1-g_{11}\,\vartheta_2}<0, \; \;u^*_1, \; u^*_2, \; \mbox{ and  } \; u^*_3, 
\]
where $u^*_1$, $u^*_2$, and $u^*_3$ are the three roots of $p(u)=0$. We shall check that $u^*_1$, $u^*_2$, and $u^*_3$ indeed are three distinct real roots using  Fan{'s}  method. As in \cite{Shengjin89Cubic}, we define by
\begin{equation*}
A:=k_2^2-3\,k_1\,k_3, \quad B:=k_1\,k_2-9\,k_0\,k_3, \quad C:=k_1^2-3\,k_0\,k_2
\end{equation*}
and the \textit{discriminant}
\begin{equation*}
\Delta_{dis}:=B^2-4\,A\,C.
\end{equation*}


\begin{lemma}[\textbf{Fan's method} \cite{Shengjin89Cubic}]\thlabel{lem: discriminant}
There are three {possible} cases using the discriminant $\Delta_{dis}$:
\begin{itemize}
  \item [(i)]   If $\Delta_{dis}>0$, then $p(u)=0$ has one real root and two nonreal complex conjugate roots.
  \item [(ii)]   If $\Delta_{dis}=0$, then $p(u)=0$ has three real roots with one root which is at least of multiplicity 2.
  \item [(iii)]  If $\Delta_{dis}<0$, then $p(u)=0$ has three distinct real roots.
\end{itemize}

\end{lemma}
We shall show that when \textbf{(H1)} is violated, then $\Delta_{dis}< 0$.  Indeed, the \texttt{Symbolic Math Toolbox} of \texttt{MATLAB} allows us to factorize $\Delta_{dis}=D_1\, D_2$, where
\[
D_1:=-3\,d_1^3\,d_2\,g_{12}^4\,g_{21}^2\,(g_{12}\,g_{21}-g_{11}\,g_{22})^2 < 0,
\]
and
\begin{align*}
D_2&:= 4\,d_1^2\,g_{12}\,g_{21}^2\,g_{22} - 27\,d_1d_2\,g_{11}^2\,g_{22}^2 + 18\,d_1d_2\,g_{11}\,g_{12}\,g_{21}\,g_{22} \\
& \quad \quad \quad + d_1d_2\,g_{12}^2\,g_{21}^2 + 4\,d_2^2\,g_{11}\,g_{12}^2\,g_{21}\\
& \;= d_1d_2 \left(g_{12}^2 \,g_{21}^2 - g_{11}^2\,g_{22}^2 \right) + 18\,d_1\,d_2\,g_{11}\,g_{22}\left(g_{12}\,g_{21}-g_{11}\,g_{22}\right)\\
& \quad \quad \quad +4\left(d_1^2\,g_{12}\,g_{21}^2\,g_{22} -2\,d_1\,d_2\,g_{11}^{\frac{1}{2}}\,g_{12}^{\frac{3}{2}}\,g_{21}^{\frac{3}{2}}\,g_{22}^{\frac{1}{2}}+d_2^2\,g_{11}\,g_{12}^2\,g_{21}\right)\\
& \quad \quad \quad +8\,d_1d_2\,g_{11}^{\frac{1}{2}}\,g_{22}^{\frac{1}{2}}\left(g_{12}^{\frac{3}{2}}\,g_{21}^{\frac{3}{2}} - g_{11}^{\frac{3}{2}}\,g_{22}^{\frac{3}{2}}\right)>0
\end{align*}
when \textbf{(H1)} fails. Thus, we can apply \thref{lem: discriminant} to confirm  that the cubic equation
\begin{equation*}
p(u)=k_3\,u^3+k_2\,u^2+k_1\,u+k_0=0
\end{equation*}
has three distinct real roots $u^*_3<u^*_2<u^*_1$. This implies that the derivative of $p(u)$ must have two distinct real roots. Due to the fact that $p(\pm\infty)=\pm\infty$ and $k_0,k_3>0$, it is easy to see that either $u^*_3<u^*_2<u^*_1<0$ or $u^*_3<0<u^*_2<u^*_1$. However, we can eliminate the case $u^*_3<u^*_2<u^*_1<0$, since ${\sigma=} \sigma(u)$ {is defined} for $u>u^{\ast}>0$. For the case $u^*_3<0<u^*_2<u^*_1$, $u^*_3$ cannot be a {critical} point of $\sigma(u)$ because $u^*_3<0$. Accordingly, there are \textit{at most two} {critical} points $u^*_1$ and $u^*_2$. We have by \eqref{eqn: asym behav of sigma(u) at infty} the asymptotic behavior $\sigma(u)\rightarrow-\infty$ as $u\rightarrow u^{\ast}$ or $u\rightarrow\infty$, which leads to the fact that the number of critical points of $\sigma(u)=0$ belonging to the interval $(u^\ast,+\infty)$ can only be \textit{odd}. As a consequence, there is only $u^*_1$ located in the interval $(u^\ast,+\infty)$, that is $u^*_3<0<u^*_2< u^\ast <u^*_1$. Moreover, the maximum of $ \sigma(u)$ is attained at $u=u^*_1$, i.e. $\displaystyle\max_{u>u^{\ast}} \sigma(u)=\sigma(u^*_1)$ (see Figure~\ref{fig:sigma}). We have the following rule to know the number of solutions of $\sigma(u)=0$ using the sign of $\sigma(u^*_1)$:

\begin{figure}[ht]
        \begin{subfigure}
                \centering
                \includegraphics[width=0.3\linewidth]{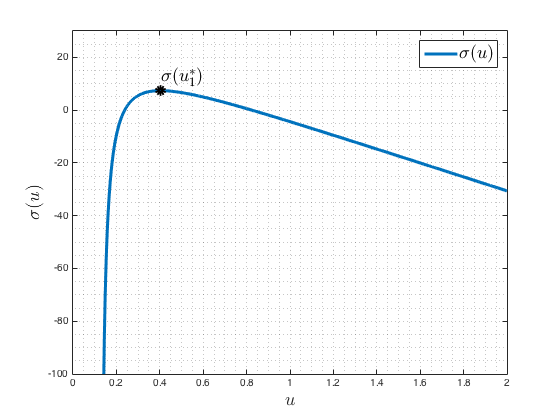}
        \end{subfigure}%
        \begin{subfigure}
                \centering
                \includegraphics[width=0.3\linewidth]{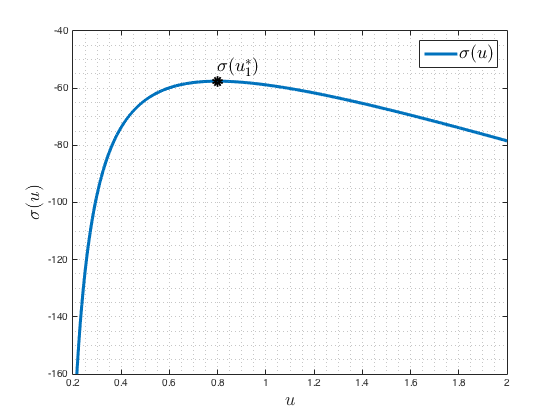}
        \end{subfigure}%
        \begin{subfigure}
                \centering
                \includegraphics[width=0.3\linewidth]{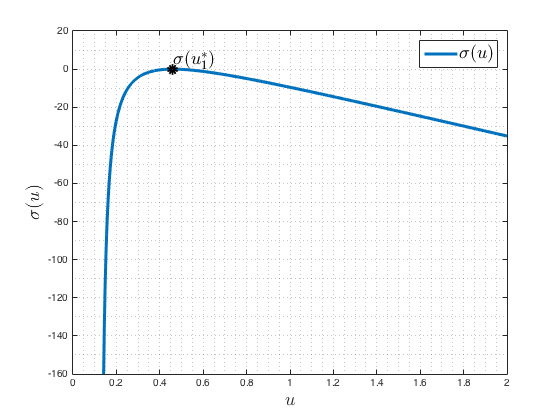}
        \end{subfigure}%
              \caption{From left to right: The graph of $\sigma(u)$ in Figure~\ref{fig: u=u(phi)} (i.e. $\sigma(u^*_1)>0$),  Figure~\ref{fig:H22} (i.e. $\sigma(u^*_1) <0$), and  Figure~\ref{fig:H23} (i.e. $\sigma(u^*_1) = 0$).}\label{fig:sigma}
\end{figure}

\begin{itemize}
  \item when $\sigma(u^*_1)>0$, equation $\sigma(u)=0$ has two distinct positive solutions;
  \item when $\sigma(u^*_1)<0$, equation $\sigma(u)=0$ has no solutions;
  \item when $\sigma(u^*_1)=0$, equation $\sigma(u)=0$ has a unique positive solution.
\end{itemize}
\begin{remark}\thlabel{re:sigma}
It follows \thref{thm: Existence of solutions to AEs} that the first two equations of \eqref{eqn: algebraic eqns for triple solns} always admit solution $(u(\phi),v(\phi))$ for all $\phi$. Hence, when we say that the system \eqref{eqn: algebraic eqns for triple solns} has no solution, we implicitly mean that $I(u(\phi),v(\phi)) \ne 0$
for all $\phi \in \mathbb{R}$. Therefore, due to the $C^1$-smoothness and positivity of $u(\phi)$ and $v(\phi)$, $I(u(\phi),v(\phi))$ must keep the same sign along the curve $(u(\phi),v(\phi))$.

Moreover, let us consider the two functions derived from \eqref{eqn: drift-diffusion DAEs algebraic eqns}, i.e. $\Phi_1(u,v):= \frac{1}{\vartheta_1}\left(c_1 - d_1 \log u - g_{11}u - g_{12}v\right),$
and $\Phi_2(u,v):= \frac{1}{\vartheta_2}\left(c_2 - d_2 \log u - g_{21}u - g_{22}v\right).$
If $(u_0,v_0,\phi_0)$  satisfies \eqref{eqn: drift-diffusion DAEs algebraic eqns}, then $\Phi_1(u_0,v_0) = \phi_0 = \Phi_2(u_0,v_0)$. Now, for $v:=v(u)$ defined in \eqref{eqn: v=v(u)}, we have
\begin{align*}
\sigma(u_0) &= \Phi_1(u_0,v(u_0)) - \Phi_2(u_0,v(u_0))\\
& = \Phi_1(u_0,v_0) + \frac{1}{\vartheta_1}g_{12}(v_0 - v(u_0)) - \Phi_2(u_0,v(u_0))\\
&= \Phi_2(u_0,v_0) + \frac{1}{\vartheta_1}g_{12}(v_0 - v(u_0)) - \Phi_2(u_0,v(u_0))\\
&= \frac{1}{\vartheta_2}\left[d_2(\log v(u_0) - \log v_0) + g_{22} (v(u_0)-v_0) \right]+ \frac{1}{\vartheta_1}g_{12}(v_0 - v(u_0)).
\end{align*}
Hence, for the case $\sigma(u^*_1) <0$, we have $\sigma(u_0)<0$ and $v_0 < v(u_0)$. This means that $I(u_0,v_0) >0$.
\end{remark}


We arrive at the following theorem.

\begin{theorem}[\textbf{Bifurcation when \textbf{(H1)} is violated}]\thlabel{thm: g11g22<g12g21 existence of alge soln}
Assume that $\textbf{(H1)}$ fails and let $(u(\phi),v(\phi))$ be a pair of solutions to \eqref{eqn: drift-diffusion DAEs algebraic eqns}. Then 
\begin{itemize}
  \item [(i)] \textbf{(triple piecewise $C^1$ solutions (cf. Figure~\ref{fig: u=u(phi)}))} when $\sigma(u^*_1)>0$, there exist $\underaccent\bar{\phi},\bar{\phi}\in\mathbb{R}$ such that 
  \begin{itemize}
  \item[$\triangleright$] for $\phi\in(-\infty,\underaccent\bar{\phi})\cup(\bar{\phi},\infty)$: $u(\phi)$ and $v(\phi)$ can be represented uniquely; are of class $C^1$; and $u'(\phi)<0$ and $v'(\phi)>0$;
  \item[$\triangleright$] at $\phi=\underaccent\bar{\phi},\bar{\phi}$: $u(\phi)$ (and $v(\phi)$) takes two distinct values;
  \item[$\triangleright$] for $\phi\in(\underaccent\bar{\phi},\bar{\phi})$: $u(\phi)$ (and $v(\phi)$) takes three distinct values $u_j(\phi)$ (and  $v_j(\phi)$), $j=1,2,3$. For each  $j$, the curve $u_j(\phi)$ (and $v_j(\phi)$) is of class $C^1$;
 \end{itemize}
  \item  [(ii)] \textbf{(unique  $C^1$-smooth monotone solutions (cf. Figure~\ref{fig:H22}))} when $\sigma(u^*_1)<0$, $u(\phi)$ and $v(\phi)$ can be represented uniquely for $\phi\in\mathbb{R}$. Moreover, $u(\phi)$ and $v(\phi)$ are of class $C^1$ with $u'(\phi)<0$ and $v'(\phi)>0$ for $\phi\in\mathbb{R}$;
   \item [(iii)] \textbf{(unique piecewise $C^1$-smooth monotone solutions (cf. Figure~\ref{fig:H23}))} when $\sigma(u^*_1)=0$, there exists $\check{\phi}\in\mathbb{R}$  such that 
   \begin{itemize}
   \item[$\triangleright$] for $\phi\in\mathbb{R}\setminus \{\check{\phi}\}$: $u(\phi)$ and $v(\phi)$ can be represented uniquely; are of class $C^1$; and  $u'(\phi)<0$ and $v'(\phi)>0$; 
    \item[$\triangleright$] at $\phi=\check{\phi}$: $u'(\phi)= -\infty$ and $v'(\phi)=\infty$.
    \end{itemize}
   \end{itemize}
\end{theorem}

\begin{proof} \textbf{Step 1:} We first check (ii). For any $\phi \in \mathbb{R}$, \thref{thm: Existence of solutions to AEs} guarantees the existence of the solution $(u(\phi),v(\phi))$ to \eqref{eqn: drift-diffusion DAEs algebraic eqns}. Following the argument in \thref{re:sigma}, we have $I(u,v)>0$
for all pair $(u,v)$ satisfying \eqref{eqn: drift-diffusion DAEs algebraic eqns}.

We now prove the uniqueness of \eqref{eqn: drift-diffusion DAEs algebraic eqns} when \textbf{(H1)} is violated for the case $\sigma(u^*_1)<0$  by contradiction. Indeed, assume that there exists $\phi_0 \in \mathbb{R}$ such that \eqref{eqn: drift-diffusion DAEs algebraic eqns} admits at least two distinct solutions $(u_1,v_1)$ and $(u_2,v_2)$. Let $M_1(u,v)$ and $M_2(u,v)$ be the functions in  \thref{thm: Uniqueness of solutions to AEs}, we have $M_1(u_1,v_1) > M_2(u_1,v_1)$ and $M_1(u_2,v_2)>M_2(u_2,v_2)$. Repeating the argument in (S3) of  \thref{thm: Uniqueness of solutions to AEs}, we get a pair $(u_3,v_3)$ satisfying 
\eqref{eqn: drift-diffusion DAEs algebraic eqns} and $M_1(u_3,v_3) < M_2(u_3,v_3)$, which is a contradiction since $I(u_3,v_3)>0$.

The $C^1$-smoothness of $u(\phi)$ and $v(\phi)$ is due to the Implicit Function Theorem, whilst the fact that $u'(\phi) <0$ and $v'(\phi)>0$ follow immediately \eqref{eqn: u'(phi)}, \eqref{eqn: v'(phi)} and the fact that $I(u(\phi),v(\phi))>0$ for all $\phi \in \mathbb{R}$.



\textbf{Step 2:} Let us consider the case (iii). Assume that $(\check{u},\check{ v}, \check{\phi})$ is the unique solution to \eqref{eqn: algebraic eqns for triple solns}. Since $\sigma(u) < 0 $ for all $u \in (u^*,\infty) \setminus\{\check{u}\}$, we can repeat the argument in (ii) for $\phi \in (-\infty,\check{\phi})\cup (\check{\phi},\infty)$ to get the assertion in (iii).

\textbf{Step 3:} For the case (i), let $(\underaccent\bar{u},\underaccent\bar{v},\underaccent\bar{\phi})$ and $(\bar{u},\bar{v},\bar{\phi})$ be two distinct solutions to \eqref{eqn: algebraic eqns for triple solns}. 
\begin{itemize}
\item For $\phi \in (-\infty,\underaccent\bar{\phi}) \cup (\bar{\phi},\infty)$, we see that $\sigma(u) <0 $ for all $u \in (u^*,\infty)\setminus (\underaccent\bar{u},\bar{u})$, in the same manner as (ii), we get the result.
\item At $\phi = \underaccent\bar{\phi}$, let $(u(\phi),v(\phi))$ be the unique solution to \eqref{eqn: drift-diffusion DAEs algebraic eqns} for $\phi \in (-\infty,\underaccent\bar{\phi})$. Let $u_1:=u(\underaccent\bar{\phi})$ and $v_1:=v(\underaccent\bar{\phi})$. We will show that $(u_1,v_1)\ne (\underaccent\bar{u},\underaccent\bar{v})$ by contradiction. Assume that $u_1=\underaccent\bar{u}$, then \eqref{eqn: drift-diffusion DAEs algebraic eqns} implies $v_1=\underaccent\bar{v}$. Thus, $(u_1,v_1)$ is the unique intersection of the two curves defined \eqref{eqn: drift-diffusion DAEs algebraic eqns} when $\phi = \underaccent\bar{\phi}$. Indeed, if there is another pair $(u_0,v_0)$ satisfying  \eqref{eqn: drift-diffusion DAEs algebraic eqns}, since $(u_0,v_0,\underaccent\bar{\phi})$ is not a solution of  \eqref{eqn: algebraic eqns for triple solns}, we have $I(u_0,v_0) \ne 0$.
Applying the Implicit Function Theorem at $\underaccent\bar{\phi}$, we get a contradiction to the uniqueness of   \eqref{eqn: drift-diffusion DAEs algebraic eqns} on $(-\infty,\underaccent\bar{\phi})$. 

Notice that $\sigma(u_1)=0$ and $\sigma'(u_1)\ne 0$ (see Figure~\ref{fig:sigma}). Taking into account \thref{re:sigma} and \eqref{eqn: v=v(u)},  since $(u_1,v_1,\underaccent\bar{\phi})$ satisfies \eqref{eqn: algebraic eqns for triple solns}, we get $\Phi_1(u_1,v_1)=\Phi_2(u_1,v_1)$ and $v_1 = v(u_1)$. Thus, in this case, $v'(u_1)= -\frac{g_{21}}{\frac{d_2}{v_1}+g_{22}}$, therefore
\begin{align*}
\sigma'(u_1) = -\frac{1}{\vartheta_1}\left(\frac{d_1}{u_1} + g_{11} +g_{12}v'(u_1)\right) + \frac{v'(u_1)}{\vartheta_2}\left(\frac{d_2}{v(u_1)} + \frac{g_{21}}{v'(u_1)} + g_{22} \right) =0,
\end{align*}
which is a contradiction. 
\item In the same manner, $u(\phi)$ (and $v(\phi)$) takes two distinct values at $\phi=\bar{\phi}$.
\item For $\phi \in (\underaccent\bar{\phi},\bar{\phi})$. Let $\underaccent\bar{u}_1$ and $\bar{u}_1$ be two distinct value of $u(\phi)$ at $\underaccent\bar{\phi}$, and let $\underaccent\bar{u}_2$ and $\bar{u}_2$ be two distinct value of $u(\phi)$ at $\bar{\phi}$. Without loss of generality, we can assume that $\sigma(\bar{u}_1)\ne 0$ and $\sigma(\underaccent\bar{u}_1)=0$. The Implicit Function Theorem at $(\underaccent\bar{\phi},\bar{u}_1)$ yields a unique $C^1$ curve $u_1(\phi)$ passing $(\underaccent\bar{\phi},\bar{u}_1)$ and satisfying \eqref{eqn: drift-diffusion DAEs algebraic eqns} for $\phi \in (\underaccent\bar{\phi},\bar{\phi})$. In $\phi$-$u$ plane, this curve $u_1(\phi)$ cuts the vertical line $\phi=\bar{\phi}$ at one of the two points $\{\underaccent\bar{u}_2,\bar{u}_2\}$. Without loss of generality, we call the intersection point $\bar{u}_2$. Since the curve $u_1(\phi)$ for $\phi \in (\underaccent\bar{\phi},\bar{\phi})$ is indeed a continuation of the unique $C^1$ curve $u(\phi)$ for $\phi \in (-\infty,\underaccent\bar{\phi})$, we have $u_1'(\phi) <0$ for $\phi \in (\underaccent\bar{\phi},\bar{\phi})$ and therefore, $\bar{u}_1 > \bar{u}_2$.

Applying \thref{prop: alge eqns} for  $\bar{u}_1$ and $\bar{u}_2$, we see that the other two points $\underaccent\bar{u}_1$ and $\underaccent\bar{u}_2$ must be simultaneously either smaller than $\bar{u}_2$ or larger than $\bar{u}_1$. If $\underaccent\bar{u}_2 > \bar{u}_1$, we can utilize  \thref{prop: alge eqns} to obtain a contradiction to the Implicit Function Theorem at the point $(\underaccent\bar{\phi},\bar{u}_1)$. Thus, $\underaccent\bar{u}_2< \bar{u}_2$.

Applying \thref{prop: alge eqns} for $\underaccent\bar{u}_2$ and $\bar{u}_2$, we must have $\underaccent\bar{u}_2< \underaccent\bar{u}_1< \bar{u}_2$. Moreover, for $\phi \in (\underaccent\bar{\phi},\bar{\phi})$, there is a $C^1$ curve $u_2(\phi)$ connecting $\bar{u}_2$ and $\underaccent\bar{u}_1$ and satisfying \eqref{eqn: drift-diffusion DAEs algebraic eqns};  and a $C^1$ curve $u_3(\phi)$ connecting $\underaccent\bar{u}_1$ and $\underaccent\bar{u}_2$ and satisfying \eqref{eqn: drift-diffusion DAEs algebraic eqns} (see Figure~\ref{fig:H21-zoom}). Employing the Implicit Function Theorem at $\bar{u}_2$, we get $\sigma(\bar{u}_2)=0$. As a consequence, $\sigma(\underaccent\bar{u}_2)\ne 0$.
\end{itemize}
\end{proof}

\begin{remark} Let $u_j(\phi)$ and $v_j(\phi)$ ($j=1,2,3)$ be the curved introduced in the proof of \thref{thm: g11g22<g12g21 existence of alge soln}. Then for each $j$, the pair $(u_j(\phi),v_j(\phi))$ solves \eqref{eqn: drift-diffusion DAEs algebraic eqns}.

\thref{thm: g11g22<g12g21 existence of alge soln} also inspires us a simple criterion to check the bifurcation of \eqref{eqn: drift-diffusion DAEs algebraic eqns}. Indeed, for any given parameters $d_1,d_2,g_{11},g_{12},g_{21},g_{22}$, we can solve the cubic equation $p(u)=0$ to get the maximum root $u^*_1$. By considering the sign of $\sigma(u^*_1)$ and taking into account  \thref{thm: g11g22<g12g21 existence of alge soln}, we can decide whether the system \eqref{eqn: drift-diffusion DAEs algebraic eqns} admits either unique $C^1$, or unique piecewise $C^1$, or triple piecewise $C^1$ solutions. 
\end{remark}


\begin{figure}[ht]
\center
\scalebox{0.55}[0.55]{\includegraphics{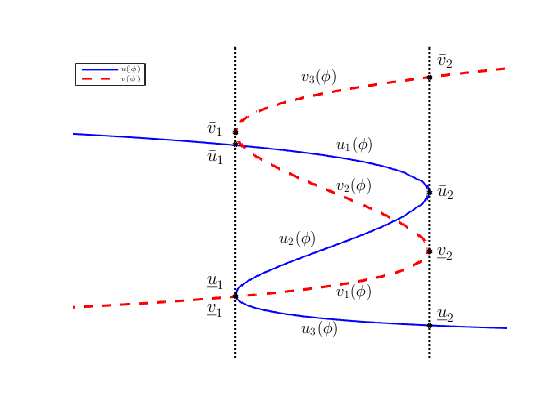}}
\caption{\textbf{A zoom-in of Figure~\ref{fig: u=u(phi)}.}
}
\label{fig:H21-zoom}
\end{figure}

\section{Auxiliary results} \label{sec:appendix}

For the reader's convenience, we quote here the definition of pseudomonotone operator in \cite{looker2006semilinear}.
\begin{definition}\label{def:pseudomonotone}
Let $A:X\to X^*$ be an operator on the real reflexive Banach space $X$. Then $A$ is pseudomonotone if and only if $u_j \rightharpoonup u$ weakly in $X$ and 
$\limsup_j \left< Au_j,u_j-u\right> \le 0,$ 
implies $\left<Au,u-v\right> \le \liminf_j \left< Au_j,u_j-v\right>$,  for all $v \in X$.
\end{definition}

We now modify the proof in \cite{looker2006semilinear} to get the following lemma.
\begin{lemma}\thlabel{le:Sn-pseudo}
Assume that $G:\mathbb{R}\to \mathbb{R}$ satisfies $\textbf{(G1)}$ and \textbf{(G2)}. For any $n > |G(0)|$, let $G_n$ as in \eqref{def:Gn}. Then the operator $S_n: H^1(\Omega) \to \left(H^1(\Omega)\right)^*$ defined by
\[
\left< S_n u,v\right>:= \int_{\Omega} G_n(u)v, \quad u,v \in H^1(\Omega),
\]
is pseudomonotone. Moreover,
\[
\left<S_nu,u\right> \ge \int_\Omega G(0) u \quad \mbox{ for all } \quad u \in H^1(\Omega).
\]
\end{lemma}
\begin{proof}
First, we shall check that $S_n$ is well-defined for all fixed $n \in \mathbb{N}$. Indeed, for all $u \in H^1(\Omega)$, the map $v \mapsto \int_{\Omega} G_n(u)v$ is linear bounded, since
\[
\left| \int_{\Omega} G_n(u)v\right| \le n \int_{\Omega}| v|  \le n|\Omega|^{\frac{1}{2}} \|v\|_{L^2(\Omega)}. 
\]

For any $u \in H^1(\Omega)$ and for all $n \in \mathbb{N}$, the fact that $G_n$ is monotone increasing leads to
\[
\left<S_nu,u\right> = \int_\Omega G_n(u)u = \int_\Omega (G_n(u)-G_n(0))u + \int_{\Omega} G_n(0)u \ge \int_{\Omega} G_n(0)u.
\]

Let $u_j$ be an arbitrary sequence in $H^1(\Omega)$ weakly converging to $u \in H^1(\Omega)$ such that $\limsup_{j} \left< S_n u_j,u_j-u\right> \le 0$. We shall check that
\[
\left< S_nu,u-v\right> \le \liminf_j \left<S_nu_j,u_j-v\right> \quad \mbox{ for all } \quad v\in H^1(\Omega).
\]
Since $u_j$ weakly converges to $u$ in $H^1(\Omega)$, the uniformly boundedness of the sequence $u_j$ in $H^1(\Omega)$ leads to
\[
\left|\left< S_n u_j,u_j - v \right> \right| \le n \int_{\Omega} |u_j-v| \le n |\Omega|^{\frac{1}{2}} \left(\|u_j\|_{L^2(\Omega)} + \|v\|_{L^2(\Omega)}\right) < \infty.
\]
Thus, up to subsequence (we still denote the subsequence by $u_j$), we get
\[
\left< S_n u_j,u_j - v \right> \to L \quad \mbox{as} \quad j \to \infty,
\]
where $L:=\liminf_k \left< S_n u_k,u_k - v \right>$. For $n >|G(0)|$, we have $G_n(0) = G(0)$ and
\begin{align*}\label{eqn:decompose-Sn}
\left<S_nu_j,u_j-v\right>&= \int_{\Omega} G_n(u_j)(u_j-v) \\
& =  \int_{\Omega} \left(G_n(u_j)-G(0) \right) u_j +  \int_\Omega G(0)u_j - \int_\Omega G_n(u_j)v.\nonumber
\end{align*} 

Since $u_j$ weakly converges to $u$ in $H^1(\Omega)$, the compact embedding $H^1(\Omega) \hookrightarrow L^2(\Omega)$ then implies the strong convergence (up to subsequence) of $u_j$ to $u$ in $L^2(\Omega)$. Since $|\Omega| < \infty$, we also have that $\int_\Omega u_j \to \int_\Omega u$ and that $u_j \to u $ strongly in $L^1(\Omega)$. Thus, we arrive at
\[
\int_{\Omega} G(0)u_j \to \int_\Omega G(0) u \quad \mbox{as} \quad j \to \infty.
\]

Besides, (up to subsequence) we can assume that $u_j$ pointwisely converges to $u$ a.e. $\Omega$. Due to the continuity  of $G_n$, the sequence $G_n(u_j)$ also pointwisely converges to $G_n(u)$ a.e. $\Omega$ as $j$ tends to infinity. Moreover, for all fixed $v \in H^1(\Omega), |G_n(u_j) v| \le n|v|$ and $n|v| \in L^1(\Omega)$. Hence, the Lebesgue's Dominated Convergence Theorem yields that $G_n(u)v \in L^1(\Omega)$ and 
\[
\int_\Omega G_n(u_j)v \to \int_{\Omega} G_n(u)v\quad \mbox{as} \quad j \to \infty. 
\]

On the other hand, we get from the monotonicity of $G_n$ that $(G_n(u_j)-G(0))u_j \ge 0$ for all $j$. Moreover, since $u_j$ is uniformly bounded in $H^1(\Omega)$, it holds that
\[
0 \le \int_{\Omega} (G_n(u_j)-G(0))u_j \le (|G(0)| + n)|\Omega|^{\frac{1}{2}}\|u_j\|_{L^2(\Omega)} < C, 
\]
for all $j$ and for some constant $C$ independent of $j$. Applying the Fatou's Lemma, we get
\[
\int_\Omega (G_n(u)-G(0))u \le \liminf_j \int_\Omega (G_n(u_j)-G(0))u_j.
\]

This completes the proof of \thref{le:Sn-pseudo}.
\end{proof}

The following lemma is due to \cite{webb1980boundary}.
\begin{lemma}\thlabel{le:Gn-and-G}
Assume that $G: \mathbb{R}\to \mathbb{R}$ satisfies \textbf{(G1)} and \textbf{(G2)}. For any $n > |G(0)|$, let $G_n$ as in \eqref{def:Gn}. Let $\{u_n\}$ be a sequence in $H^1(\Omega)$ weakly converging to some $u$ in $H^1(\Omega)$ and satisfying 
\[
0 \le -\int_\Omega G_n(u_n)u_n \le C \quad \mbox{ for some } C>0 \quad \mbox{ and for all } \quad n.
\] 
Then, $G(u)u \in L^1(\Omega)$ and $G_n(u_n) \to G(u)$ strongly in $L^1(\Omega)$.
\end{lemma}
\begin{proof}
We modify the proof in \cite{webb1980boundary} to obtain the desired result. 

As $u_n$ weakly converges to $u$ in $H^1(\Omega)$, the sequence $u_n$ is uniformly bounded in $H^1(\Omega)$, and (up to subsequence), we can assume that $u_n$ pointwisely converges to $u$ a.e. $\Omega$. Thus, $G_n(u_n)$ also pointwisely converges to $G(u)$ a.e. $\Omega$. 

On the other hand, the monotonicity of $G$ and the fact that $u_n$ is uniformly bounded yield
\[
0 \le \int_\Omega (G_n(u_n)-G(0))u_n \le C + |G(0)|\int_\Omega |u_n| \le C + |G(0)| |\Omega|^{\frac{1}{2}} \|u_n\|_{L^2(\Omega)} \le C.
\]
Here, the constant $C$ may change from lines to lines. Applying Fatou's Lemma, we get
\[
0 \le \int_\Omega (G(u) -G(0)) u \le   \liminf_{n \to \infty} \int_\Omega (G_n(u_n)-G(0))u_n \le C.
\]
Hence,
\[
\int_\Omega |G(u)u| \le \int_\Omega (G(u)u-G(0)u) + \int_\Omega |G(0)u| < \infty. 
\]

Now, for any $\delta >0$
\[
|G_n(u_n)-G(0)| \le \sup_{|t|\le \delta^{-1}} |G(t)| + |G(0)| + \delta(G_n(u_n)-G(0))u_n.
\]
Given $\eps >0$, we can choose $\delta>0$ such that for any $E \subset \Omega$ with $|E| < \delta$, 
\[
\int_E |G_n(u_n)| \le \left( \sup_{|t|\le \delta^{-1}} |G(t)| + 2|G(0)|\right)|E| + \delta C < \eps.
\]
Here, $|E|$ denotes the Lebesgue measure of $E$. By Vitali's Convergence Theorem, we have $G_n(u_n) \to G(u)$ strongly in $L^1(\Omega)$.
\end{proof}

The following classical result \cite{looker2006semilinear,zeidler1990nonlinear} guarantees the existence of weak solution to semilinear elliptic differential equation with pseudomonotone operator.
\begin{theorem}\thlabel{thm:existence-pseudo}
Let $A: X \to X^*$ be a pseudomonotone, bounded and coercive operator on the real, separable and reflexive Banach space $X$. Then for each $b \in X^*$, the equation $Au=b, u \in X,$ has a solution.
\end{theorem}

\subsubsection*{Acknowledgements}
The authors are grateful to the anonymous referees for many helpful comments and valuable suggestions on this paper. L.-C. Hung would like to thank Professors Tai-Chia Lin and Chun Liu for introducing the problem to him. He is also grateful for their fruitful discussions and valuable comments in preparation of the manuscript and for suggesting improvements. The authors also thanks Professor Robert Eisenberg for introducing them the biological aspect of the ion channel problem and for his interest in this work. The research of L.-C. Hung is partly supported by the grant 106-2115-M-011-001-MY2 of Ministry of Science and Technology, Taiwan.

 







\end{document}